\documentclass[12pt
]{amsart}
\usepackage{amssymb}
\usepackage[dvips]{graphics}
\ifx\mathrm\undefined\let\mathrm\rm\fi
\ifx\mathbf\undefined\let\mathbf\bf\fi
\ifx\mathfrak\undefined\let\mathfrak\frak\fi
\ifx\mathcal\undefined\let\mathcal\cal\fi
\ifx\mathbb\undefined\let\mathbb\Bbb\fi
\ifx\emph\undefined\let\emph\it\fi
\font\bb=msbm10 at9.98pt
\def\semidirect{\hbox{$\;$\bb\char'156$\;$}}
\newcommand{\SL}{\mathrm{SL}}

\newcommand{\Z}{{\mathbb Z}}

\newcommand{\R}{{\mathbb R}}
\newcommand{\C}{{\mathbb C}}
\newcommand{\Q}{{\mathbb Q}}

\newcommand{\Ref}[1]{{(\ref{#1})}}

\newcommand{\LLambda}{\mathbf{\Lambda}}
\newcommand{\lamdba}{\lambda}
\newcommand{\zz}{{\vec{z}}}
\newcommand{\F}{\mathcal{F}}
\newcommand{\sig}{{p}}

\newcommand{\dontprint}[1]
{\relax}

\newtheorem%
{thm}{Theorem}[section]
\newtheorem%
{proposition}[thm]{Proposition}
\newtheorem%
{lemma}[thm]{Lemma}
\newtheorem%
{lemmadef}[thm]{Lemma-Definition}
\newtheorem%
{corollary}[thm]{Corollary}
\newtheorem%
{conjecture}[thm]{Conjecture}

\begin{document}
\title[{}]{q-deformed KZB heat equation:
completeness, modular properties and SL(3,$\Z$)}
\author[{}]
{Giovanni Felder${}^{*,1}$ 
\and Alexander Varchenko${}^{**,2}$}
\thanks{${}^1$Supported in part by the Swiss National Science Foundation} 
\thanks{${}^2$Supported in part by NSF grant  DMS-9801582}
\maketitle 
\medskip \centerline{\it ${}^*$Department of Mathematics,
  ETH-Zentrum,} \centerline{\it 8092 Z\"urich, Switzerland} \medskip
\centerline{\it ${}^{**}$Department of Mathematics, University of
  North Carolina at Chapel Hill,} \centerline{\it Chapel Hill, NC
  27599-3250, USA} \medskip 

\centerline{October 2001}
\begin{abstract}
We study the properties of
one-dimensional hypergeometric integral solutions of the
$q$-difference (``quantum'') analogue of the 
Knizhnik--Zamolodchikov--Bernard
equations on tori. We show that they also obey a difference
KZB heat equation in the modular parameter, give formulae
for modular transformations, and prove a completeness result,
by showing that the associated Fourier transform is invertible.
These results are based on $\SL(3,\Z)$ transformation properties
parallel to those of elliptic gamma functions.
\end{abstract}

\dontprint{
We study the properties of
one-dimensional hypergeometric integral solutions of the
q-difference ("quantum") analogue of the 
Knizhnik-Zamolodchikov-Bernard
equations on tori. We show that they also obey a difference
KZB heat equation in the modular parameter, give formulae
for modular transformations, and prove a completeness result,
by showing that the associated Fourier transform is invertible.
These results are based on SL(3,Z) transformation properties
parallel to those of elliptic gamma functions.
}

\section{Introduction}
In this paper we continue the study of the q-analogue of the
Knizhnik--Zamolodchikov--Bernard (qKZB) equations on elliptic curves
and their solutions initiated in \cite{FTV1,FTV2,FV1}.

In \cite{FTV1}, hypergeometric solutions of qKZB equations were
introduced.  In \cite{FTV2}, the monodromy of hypergeometric solutions
was calculated, and a symmetry between equations and monodromies was
discovered: the equations giving the monodromy are again qKZB
equations but with modular parameter and step of the difference
equations exchanged.  In \cite{FV1} we introduced the $q$-analogue of
the KZB heat equation, which governs the change in the modular
parameter of the elliptic curve, proved that it is compatible with the
other equations and explained how to recover the KZB heat equation in
the semiclassical limit.

In this paper we prove three results about our
hypergeometric solutions in the case where the sum of 
highest weights is two: we show that the hypergeometric solutions
also obeys the qKZB heat equation of \cite{FV1}, see Theorem
\ref{t-HE1}.  We give a formula (Theorem \ref{t-mod}) for the
transformation properties of the hypergeometric solutions under the
modular group $\SL(2,\Z)$. We prove a completeness result, Corollary
\ref{completeness}, by showing that the associated ``Fourier
transform'' is invertible.

Then we show that these results are part of
a bigger picture, in which the modular group combines with the
transformation defined by the qKZB heat equation to give a
set of quadratic identities for our generalized hypergeometric integrals.
In fact, this picture can already be seen in a simpler situation, in which
the elliptic gamma function \cite{R}, \cite{FV2} plays the role
of the hypergeometric integral. The elliptic gamma function
is a function $\Gamma(z,\tau,\sig)$ of three complex variables
obeying identities \cite{FV2} involving its values at points related
by an action of $\SL(3,\Z)\semidirect\Z^3$.
These identities mean that $\Gamma$ is a ``degree 1'' generalization
of a Jacobi modular function, see \cite{FV2}.

These identities are a scalar version of the identities obeyed
by the hypergeometric solutions of the qKZB equation.
Their meaning is that the hypergeometric integrals define
a discrete projectively flat
$\SL(3,\Z)$-connection (i.e., a lift of the
action to the projectivization) on a vector bundle
over ``regular'' orbits of $\SL(3,\Z)$ acting on the
variety of pairs (point in $\C^3-\{0\}$, 
plane through $0$ containing
it). This is the content of Theorem \ref{t-pf}.
The results on the elliptic gamma function are 
used in the proof, since the
``phase function'' which appears in the integrand of
hypergeometric solutions is a ratio of gamma functions. In
fact, we see ``experimentally'' that there seems to be  a principle
stating that to every identity obeyed by the gamma function,
there corresponds an identity for hypergeometric integrals.
The proofs of the identities consist in applying the gamma function
identity to the phase function in the integrand, and then use a
version of Stokes' theorem to relate the integrals. The second
step is relatively simple in the case of the one-dimensional integrals
to which we restrict ourselves here, but becomes exceedingly involved in the
higher dimensional case. Proving our identities 
in the higher dimensional case, i.e., if the sum of highest weights 
is larger than two,  remains a challenging open problem.
An alternative approach to this problem is based on representation
theory: in \cite{EV} a representation-theoretic interpretation of 
a degenerate version of the
qKZB equations was established. It was shown that traces of intertwining
operators for quantum groups satisfy 
a 
version of the qKZB equations 
and are eigenfunctions of analogues of Macdonald operators. This fact indicates
that our Theorem \ref{t-HE1} is an elliptic analogue of the Macdonald-Mehta identities
proved by Cherednik, \cite{C}, \cite{EK}, and the present work
is a glimpse 
into
an elliptic version of  Macdonald's theory.

In fact Theorem \ref{t-pf} concerns the case of generic 
parameters, where an infinite-dimensional vector bundle
is preserved by the projectively flat connection. 
In a next paper we will restrict our attention to special 
$\SL(3,\Z)$-orbits. The projectively flat connection can be
then defined on a finite-dimensional vector bundle of 
theta-functions which are a $q$-deformed version of the
space of conformal blocks in conformal field theory. In this
setting the analogy with Macdonald's theory will appear
more explicitly.
Another degeneration of the $\SL(3,\Z)$ symmetry of our hypergeometric
integrals are the $\SL(3,\Z)$ symmetries 
of the ordinary Fourier transform     
indicated in \cite{FV4}.

\section{Hypergeometric solutions of the qKZB equations}\label{s-hyp}

We use the definitions and notations of \cite{FV1}.  The elliptic
$sl_2$ qKZB equations are a compatible system of difference equations
for a function $v(\zz,\lambda,\tau)$ of $\zz\in\C^n$, $\lambda\in\C$
and $\tau$, $\mathrm{Im}\,{}\tau>0$, taking values in the zero weight
space $V_\LLambda[0]$ of a tensor product of $E_{\tau,\eta}(sl_2)$
evaluation modules.  This space comes with a basis of eigenvectors of
commuting operators $h^{(i)}$, $i=1,\ldots,n$ and depends on
parameters $\LLambda=(\Lambda_1,\dots,\Lambda_n) \in\C^n$.  The qKZB
equations have the form
\begin{equation}\label{e-qKZB}
  v(\zz+p\delta_i,\tau)=K_i(\zz,\tau,p,\eta)v(\zz,\tau), \qquad
  i=1,\dots, n.
\end{equation}
They are supplemented by the qKZB heat equation
\begin{equation}
  \label{e-heat}
  v(\zz,\tau)=T(\zz,\tau,p,\eta)v(\zz,\tau+p).
\end{equation}
The step $p$ is a complex parameter and $(\delta_i)_{i=1,\ldots,n}$ is
the standard basis of $\C^n$.  Here $v$ is viewed as a function of
$\zz$ and $\tau$ with values in the space $\mathcal{F}(V_\LLambda[0])$
of meromorphic $V_\LLambda[0]$-valued functions of a complex variable
$\lambda$.  The qKZB operators $K_i(\zz,\tau,p,\eta)$ are difference
operators in $\lambda$ and can be expressed in terms of a product of
(dynamical) $R$-matrices.  The last equation is the qKZB heat equation
and involves the integral operator $T(\zz,\tau,p,\eta)$.  The latter
is expressed in terms of hypergeometric integrals: let
$u(\zz,\lambda,\mu,\tau,p,\eta)\in V_{\LLambda}[0]\otimes
V_{\LLambda}[0]$ be the universal hypergeometric function as in
\cite{FV1}.  We may view it as a function $u(\zz,\tau,p,\eta)$ taking
its values in the $V_{\LLambda}[0]\otimes V_{\LLambda}[0]$-valued
functions of the two ``dynamical variables'' $\lambda$ and $\mu$. Then
it is a projective solution (i.e., a solution up to a constant factor)
of the qKZB equations in the first factor, and of the mirror qKZB
equations in the second:
\begin{eqnarray}\label{eq1}
  u(\zz+\delta_jp,\tau,p,\eta)&=&K_j(\zz,\tau,p,\eta)\otimes D_j\,
  u(\zz,\tau,p,\eta), \notag\\ u(\zz+\delta_j\tau,\tau,p,\eta)&=&
  D^\vee_j\otimes K^\vee_j(\zz,p,\tau,\eta)\, u(\zz,\tau,p,\eta),\\ 
  u(\zz+\delta_j,\tau,p,\eta)&=&u(\zz,\tau,p,\eta).  \notag
\end{eqnarray}
The mirror qKZB operators $K_i^\vee(\zz,p,\tau,\eta)$ are obtained
from the qKZB operators by exchanging $\tau$ and $p$ and ``reversing
the order of factors'', see \cite{FV1}.  The operators $D_j$,
$D_j^\vee$ are operators of multiplication by certain functions of the
$h^{(i)}$ and the dynamical variable $\mu$ and $\lambda$,
respectively: in terms of the function
\[
\alpha(\lambda)=\exp(-{\pi i\lambda^2/4\eta}),
\]
we have, for $j=1,\dots, n$,
\begin{eqnarray*}
  D_j(\mu)&=& \frac {\alpha(\mu-2\eta(h^{(j+1)}+\cdots+h^{(n)}))}
  {\alpha(\mu-2\eta(h^{(j)}+\cdots+h^{(n)}))} \, e^{\pi
    i\eta\Lambda_j(\sum_{l=1}^{j-1}\Lambda_l-
    \sum_{l=j+1}^{n}\Lambda_l)}, \\ D^\vee_j(\lambda)&=& \frac
  {\alpha(\lambda-2\eta(h^{(1)}+\cdots+h^{(j-1)}))}
  {\alpha(\lambda-2\eta(h^{(1)}+\cdots+h^{(j)}))} \, e^{-\pi
    i\eta\Lambda_j(\sum_{l=1}^{j-1}\Lambda_l-
    \sum_{l=j+1}^{n}\Lambda_l)}.
\end{eqnarray*}
The integral operator $T(\zz,\tau,p,\eta)$ is then
\[
T(\zz,\tau,p,\eta)v=(\alpha\otimes Q_{\tau+p})
u(\zz,\tau,\tau+p,\eta)\otimes v.
\]
$\alpha$ is the operator of multiplication by the function
$\alpha(\lambda)$ and $Q_{\tau+p}$ is a bilinear form on
$V_\LLambda[0]$-valued functions, whose kernel is the Shapovalov
bilinear form on $V_\LLambda[0]$:
\[
Q_{\sigma}(f\otimes g)=\int Q(\mu,\sigma,\eta)(f(\mu),g(-\mu))
\alpha(\mu)d\mu.
\]

\begin{figure}
  \rotatebox{-90}{\scalebox{.5}{\includegraphics{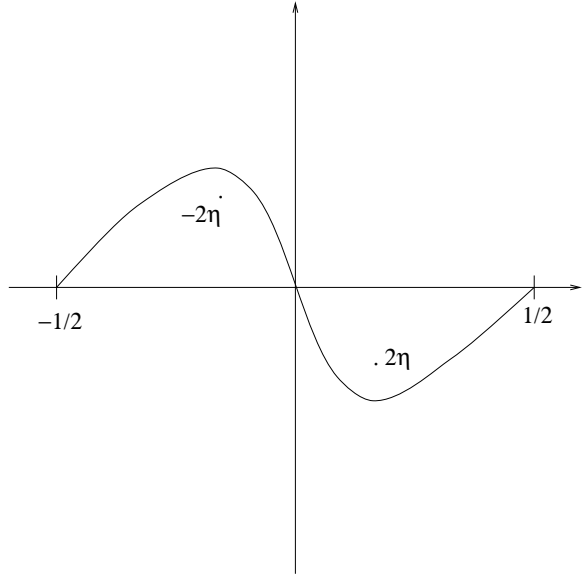}}}
\caption{The integration cycle}\label{Fig1}
\end{figure}

If $\Lambda_1+\cdots+\Lambda_n=2$, the universal hypergeometric
function is given in terms of Jacobi's first theta function $\theta$,
and the phase function, see Appendix \ref{App-A}, by the following formulae:
$V_\LLambda[0]$ has a basis $\epsilon_j=e_0\otimes\cdots\otimes
e_1\otimes\cdots\otimes e_0$, with $e_1$ in the $j$th factor
($j=1,\dots,n$) and
$h^{(i)}\epsilon_j=(\Lambda_j-2\delta_{ij})\epsilon_j$. Then the {\em
  weight functions} are $\omega(t,\zz,\lambda,\tau,\eta)
=\sum_{a=1}^n\omega_a(t,\zz,\lambda,\tau,\eta) \epsilon_a$, with
\[
\omega_a(t,\zz,\lambda,\tau,\eta)= \frac
{\theta(\lambda+2\eta+t-z_a-\eta\Lambda_a-2\eta\sum_{j<a}\Lambda_j,\tau)}
{\theta(t-z_a-\eta\Lambda_a,\tau)} \prod_{j=1}^{a-1} \frac
{\theta(t-z_j+\eta\Lambda_j,\tau)}
{\theta(t-z_j-\eta\Lambda_j,\tau)}\,.
\]
On the other hand, $\omega^\vee(t,\zz,\mu,p,\eta)=\sum
\omega^\vee_b(t,\zz,\mu,p,\eta)\epsilon_b$, where
\[
\omega^\vee_b(t,\zz,\mu,p,\eta)= \frac
{\theta(\mu+2\eta+t-z_b-\eta\Lambda_b-2\eta\sum_{j>b}\Lambda_j,p)}
{\theta(t-z_b-\eta\Lambda_b,p)} \prod_{j=b+1}^{n} \frac
{\theta(t-z_j+\eta\Lambda_j,p)} {\theta(t-z_j-\eta\Lambda_j,p)}\,.
\]
The universal hypergeometric function $u$ is expressed in terms of
these weight functions and the phase function
\[
\Omega_a(z,\tau,p)=\prod_{j,k=0}^\infty \frac {(1-e^{2\pi
    i(z-a+j\tau+kp)})(1-e^{2\pi i(-z-a+(j+1)\tau+(k+1)p)})}
{(1-e^{2\pi i(z+a+j\tau+kp)})(1-e^{2\pi i(-z+a+(j+1)\tau+(k+1)p)})}
\,.
\]
We then have $u(\zz,\lambda,\mu,\tau,p,\eta)=$
$\sum_{a,b}u_{a,b}(\zz,\lambda,\mu,\tau,p,\eta)\epsilon_a\otimes
\epsilon_b$, with
\begin{equation}\label{eq-467}
  u_{a,b}(\zz,\lambda,\mu,\tau,p,\eta) =
  e^{-\frac{i\pi\lambda\mu}{2\eta}} \int_\gamma
  \prod_{j=1}^n\Omega_{\eta\Lambda_j}(t-z_j,\tau,p)
  \omega_a(t,\zz,\lambda,\tau,\eta) \omega^\vee_b(t,\zz,\mu,p,\eta)dt.
\end{equation}
The integral is defined as the analytic continuation from the region
where Re$(\Lambda_i)<0$, $z_j\in\R$ and Im$(\eta)<0$. In this region,
the integration cycle is just the interval $[0,1]$.  After the
analytic continuation, the cycle is deformed to go above the pole at
$z_j-\eta\Lambda_j$ and below the pole at $z_j+\eta\Lambda_j$, as in
Fig.~\ref{Fig1}.

The Shapovalov form is $Q(\lambda,\tau,\eta)(\epsilon_a,
\epsilon_b)=\delta_{a,b} Q_a(\lambda,\tau,\eta)$ with
\[
Q_a(\lambda,\tau,\eta) =
\frac{\theta(2\eta\Lambda_a,\tau)\theta'(0,\tau)}
{\theta(\lambda-2\eta+2\eta\sum_{j<a}\Lambda_j,\tau)
  \theta(\lambda-2\eta+2\eta\sum_{j\leq a}\Lambda_j,\tau)}\,.
\]
Our first result is that our hypergeometric projective solutions of
the qKZB equations are also projective solutions of the quantum heat
equation:

\begin{thm}\label{t-HE1} Suppose that $\Lambda_1+\cdots+\Lambda_n=2$.
  Let us view $u(\zz,\lambda,\mu,\tau,p,\eta)$ as a function
  $u(\zz,\tau,p,\eta)$ with values in the space of
  $V_{\LLambda}[0]\otimes V_{\LLambda}[0]$-valued functions of
  $\lambda$ and $\mu$.  Then $u$ is a projective solution of the qKZB
  heat equation
\[u(\zz,\tau,p,\eta)=
C\,T(\zz,\tau,p,\eta)\otimes D_0\,u(\zz,\tau+p,p,\eta),\] where
\[C=-\frac{e^{4\pi i\eta}}
{2\pi\sqrt{4i\eta}}\,.
\]
Here $T(\zz,\tau,p,\eta)$ acts on the first dynamical variable
$\lambda$ and $D_0$ is the operator of multiplication by the function
$e^{-\pi i \mu^2/4\eta}$ of the second dynamical variable $\mu$.
\end{thm}

\medskip

It is then easy to construct true solutions to the system
\Ref{e-qKZB}, \Ref{e-heat} from these projective solutions: for any
$\mu$, $1\otimes\prod_jD_j(\mu)^{-z_j/p}D_0(\mu)^{\tau/p}
u(\zz,\lambda,\mu,\tau,p,\eta)C^{\tau/p}$, viewed as a function of
$\zz,\lambda,\tau$, obeys \Ref{e-qKZB} and \Ref{e-heat} in the first
factor.

In more explicit terms, we have the following statement.

\medskip

\begin{thm}\label{q-heat-n}  Let
  for $a=1,\dots,n$
\[
\rho_a(\lambda,\zz,\tau,\eta)=e^ { \frac{2\pi
    i\eta}\tau(\sum_{1}^{a-1}z_j\Lambda_j-\sum_{a}^nz_j\Lambda_j)
  -\frac{i\pi}\tau(\lambda+2\eta-2\eta\sum_{1}^a\Lambda_j-2z_a)
  (\lambda+2\eta-2\eta\sum_{1}^{a-1}\Lambda_j)
  -\frac{i\pi}\tau\lambda\sum_{1}^nz_j\Lambda_j },
\]
and let the {\em fundamental hypergeometric solution} $\bar u= \sum
\bar u_{a,b}\epsilon_a\otimes\epsilon_b$ be defined by
\[
\bar u_{a,b}(\zz,\lambda,\mu,\tau,p,\eta)
=e^{-\frac{i\pi\mu^2\tau}{4\eta p}}
\rho_{b}(-\mu,\zz,p,\eta)u_{a,b}(\zz,\lambda,\mu,\tau,p,\eta).
\]
Then, for any $\mu,b$, the $V_\LLambda[0]$-valued function
$v(z,\lambda,\tau)=\sum_a\bar
u_{a,b}(\zz,\lambda,\mu,\tau,p,\eta)\epsilon_a$ is a solution of the
system
\begin{eqnarray*}
  v(\zz+p\delta_i,\tau)&=&K_i(\zz,\tau,p,\eta)v(\zz,\tau), \qquad
  i=1,\dots, n, \\ v(\zz,\tau)&=&C\,T(\zz,\tau,p,\eta)v(\zz,\tau+p),
\end{eqnarray*}
where $C=-e^{4\pi i\eta}/2\pi\sqrt{4 i\eta}$.
\end{thm}

\medskip

\noindent{\it Proof:}
A consequence of \Ref{eq1}, Theorem \ref{t-HE1} and the identity
\begin{equation}\label{e-d99}
  \prod_{i=1}^n D_i(\mu)^{-\frac {z_i}{p}}\epsilon_b =
  \rho_b(-\mu,\zz,p,\eta)e^{\frac{i\pi}p
    (\mu-2\eta+2\eta\sum_1^b\Lambda_j)
    (\mu-2\eta+2\eta\sum_1^{b-1}\Lambda_j)} \epsilon_b.
\end{equation}
The exponential function on the right-hand side is independent of
$\zz,\lambda,\tau$ and therefore does not affect the statement of the
theorem.  \hfill$\square$

\medskip

\noindent{\bf Remarks.}
\begin{enumerate}
\item The constant $C$ could also be eliminated by including a factor
  $C^{\tau/p}$ in $\bar u$, but it is simpler to consider $C$ as part
  of the heat equation.
\item In the definition of the fundamental hypergeometric solution
  above, we have included an additional factor, see \Ref{e-d99}, with
  respect to the obvious choice. This leads to simpler formulae in the
  next section, since $\rho_b$ appears in the modular transformations
  of the qKZB operators.
\end{enumerate}

\section{Modular transformations}\label{s-modt}
The coefficients of the qKZB equations are quasi-periodic functions of
$z_1,\dots,z_n,\lambda$ with periods $1$ and $\tau$. It is therefore
natural to consider the qKZB equations as equations for sections of a
certain vector bundle over a Cartesian power of the elliptic curves
$\C/\Z+\tau\Z$.
 
Since, for any $\left(\begin{array}{cc}a&b\\ c&d\end{array}\right)\in
\SL(2,\Z)$, the elliptic curves with modulus $\tau$ and $(a\tau+b)/
(c\tau+d)$ are isomorphic, the corresponding qKZB equations are
related. For generators, the formulae relating solutions are the
following. Let for $a=1,\dots,n$
\[
B_a(\zz,\lambda,\tau,p,\eta)= e^{-\frac{\pi i\eta}{\tau p}
  \sum_{j<k}(z_j-z_k)^2\Lambda_j\Lambda_k +\frac{\pi i
    p}{4\eta\tau}\lambda^2}\rho_a(\lambda,\zz,\tau,\eta).
\]
Then we have:
\begin{proposition}
\ 
  \begin{enumerate}
  \item[(i)] Suppose that
    $v(\zz,\lambda)=\sum_{a=1}^nv_a(\zz,\lambda)\epsilon_a$ is a
    solution of the qKZB equations with parameters $\tau+1,p,\eta$.
    Then $v(\zz,\lambda)$ is a solution of the qKZB equations with
    parameters $\tau, p,\eta$.
  \item[(ii)] Suppose that
    $v(\zz,\lambda)=\sum_{a=1}^nv_a(\zz,\lambda)\epsilon_a$ is a
    solution of the qKZB equations with parameters
    $-1/\tau,p/\tau,\eta/\tau$. Then
   \[\tilde v(\zz,\lambda)=\sum_{a=1}^nB_a(\zz,\lambda,\tau,p,\eta)
   v_a(\zz/\tau,\lambda/\tau) \epsilon_a\] is a solution of the qKZB
   equations with parameters $\tau, p,\eta$.
  \end{enumerate}
\end{proposition}

The proof of this proposition is based on the formulae for
transformation properties of the qKZB equations under modular
transformations, see Appendix \ref{App-B}.

Now the question about monodromy is well-posed: can one express the
fundamental solution at the transformed values of the parameters in
terms of the fundamental solution at the original values? The answer
is provided by the following result.

\begin{thm}
\ 
  \begin{enumerate}
  \item[(i)]
$\bar u_{a,b}(\zz,\lambda,\mu,\tau+1,p,\eta)=
e^{-\frac{i\pi\mu^2}{4\eta p}}\bar
u_{a,b}(\zz,\lambda,\mu,\tau,p,\eta).
$
\item[(ii)] Suppose that $\mathrm{Im}(\eta\tau/p)<0$,
  $\mathrm{Im}(p/\tau)>0$.
  \[
  B_a(\zz,\lambda,\tau,p,\eta) \bar u_{a,b}\left(\textstyle{ \frac
      \zz\tau, \frac\lambda\tau, \frac\nu\tau, -\frac1\tau, \frac
      p\tau, \frac\eta\tau }\right) = \sum_{c=1}^n\int \bar
  u_{a,c}(\zz,\lambda,\mu,\tau,p,\eta)
  M_{c,b}(\zz,\mu,\nu,\tau,p,\eta)d\mu.
\]
The monodromy matrix $(M_{c,b})$ is
\begin{eqnarray*}
  M_{c,b}(\zz,\mu,\nu,\tau,p)&=& e^{\frac{2\pi i\eta}{3\tau p}
    (\eta^2\sum_1^n\Lambda_j^3+\tau^2+p^2-3p+3\tau+3\tau p+1)
    -\frac{2\pi i}{p\tau} (\nu-2\eta+2\eta\sum_1^{b-1}\Lambda_j)
    (\nu-2\eta+2\eta\sum_1^{b}\Lambda_j)} \\ &&
  \times\frac1{2\pi i}\sqrt{\frac{ip}{4\eta\tau}} Q_c(\mu,p,\eta)u_{c,b}
  \left(\textstyle{ \frac \zz p, -\frac\mu p, \frac \nu p, -\frac1p,
      -\frac\tau p, \frac\eta p}  \right).\end{eqnarray*} The
integration over $\mu$ is over the path $x\mapsto x\eta+\epsilon$,
$x\in\R$ for any generic real $\epsilon$.
  \end{enumerate}
\end{thm}

\medskip

The first statement of the  theorem 
is trivial to check. The second is the first
of the two {\em modular relations}.
Introduce functions $\rho, \rho^{\vee}$:
\begin{eqnarray*}
  \rho_a(\lambda,\zz,\tau,\eta) &=&e^ { \frac{2\pi
      i\eta}\tau(\sum_{1}^{a-1}z_j\Lambda_j-\sum_{a}^nz_j\Lambda_j)
    -\frac{i\pi}\tau(\lambda+2\eta-2\eta\sum_{1}^a\Lambda_j-2z_a)
    (\lambda+2\eta-2\eta\sum_{1}^{a-1}\Lambda_j)
    -\frac{i\pi}\tau\lambda\sum_{1}^nz_j\Lambda_j },\\ 
  \rho_b^\vee(\lambda,\zz,\tau,\eta) &=&e^ { \frac{2\pi
      i\eta}\tau(\sum_{b+1}^{n}z_j\Lambda_j-\sum_{1}^bz_j\Lambda_j)
    -\frac{i\pi}\tau(\lambda+2\eta-2\eta\sum_{b}^n\Lambda_j-2z_b)
    (\lambda+2\eta-2\eta\sum_{b+1}^{n}\Lambda_j)
    -\frac{i\pi}\tau\lambda\sum_{1}^nz_j\Lambda_j }.
\end{eqnarray*}

Then the modular relations are the identities:
\begin{thm}\label{t-mod}\ 

\noindent{\rm (i)}
 Suppose that $\mathrm{Im}(\eta \tau/p) < 0$,
  $\mathrm{Im}(p/\tau)>0$. Then the universal hypergeometric function
$u$ satisfies the following relation,

\begin{eqnarray}\label{eq-claim5}
  \sum_{c=1}^n\int
\lefteqn{u_{a,c}
    (\zz,\lambda,\mu,\tau,p,\eta) u_{c,b} \left(\textstyle{ \frac \zz
        p, -\frac\mu p, \frac \nu p, -\frac1p, -\frac\tau p, \frac\eta
        p} \right)}\nonumber \\ & &\times
  Q_c(\mu,p,\eta)\rho_c(-\mu,\zz,p,\eta) e^{-\frac{\pi
      i\tau\mu^2}{4\eta p}}\,d\mu\nonumber \\ 
  &=&2\pi i\sqrt{\frac{4\eta\tau}{ip}} \rho_a(\lambda,\zz,\tau,\eta)
  \rho_b^\vee \left(\textstyle{ \frac \nu p, \frac \zz p, -\frac\tau
      p, \frac\eta p}\right) e^{\frac{i\pi
      p}{4\eta\tau}(\lambda^2+(\nu/p)^2)} \\ & &\times
  u_{a,b}\left(\textstyle{ \frac \zz\tau, \frac\lambda\tau,
      \frac\nu\tau, -\frac1\tau, \frac p\tau, \frac\eta\tau
      }\right)e^{-\frac{\pi i\eta}{3p\tau}\psi}, \nonumber \\ \psi&=&
  3\sum_{j<k}\Lambda_j\Lambda_k(z_j-z_k)^2
  +2\left(\sum_{j=1}^n\eta^2\Lambda_j^3 +\tau^2+p^2-3p+3\tau+3\tau p
    +1\right).\nonumber
\end{eqnarray}
The integration over $\mu$ is over the path $x\mapsto x\eta+\epsilon$,
$x\in\R$ for any generic real $\epsilon$.

\noindent{\rm(ii)} Suppose that $\mathrm{Im}(\eta p/\tau) < 0$,
  $\mathrm{Im}(\tau/p)>0$. Then the universal hypergeometric function
$u$ satisfies the following relation,
$$
 \sum_{c=1}^n\int  
u_{a,c}\left(\textstyle{ \frac 
      \zz\tau, \frac\lambda\tau, \frac\mu\tau, -\frac p\tau, 
-\frac1\tau, \frac\eta\tau }\right)\, 
u_{c,b}(\zz,-\mu,\nu,\tau,p,\eta)
\,Q_c(\mu,\tau,\eta) \,
\rho^\vee_c (\mu,\zz,\tau,\eta)\,
e^{-{i\pi p\over 4\eta \tau}\mu^2}\,d\mu \,=
$$
$$ 
2\pi i 
\sqrt{\frac{4\eta p}{i\tau}}\, 
\rho_a \left(\textstyle{\frac\lambda p,\frac\zz p,\frac\tau p,
\frac\eta p} \right)^{-1}\,
  \rho_b^\vee (  \nu,  \zz,  p, \eta ) \,
u_{a,b}\left(\textstyle{ \frac\zz p, 
\frac\lambda p, \frac\nu p, \frac \tau p, 
-\frac1 p , \frac\eta p}\right)\, 
e^{{i\pi \tau \over 4\eta  p}((\lambda/\tau)^2 + \nu^2)}
\, e^{-{\pi i\eta \over 3p\tau}\psi}, 
$$
$$
 \psi  = 
  3\sum_{j<k}\Lambda_j\Lambda_k(z_j-z_k)^2
  +2\left(\sum_{j=1}^n\eta^2\Lambda_j^3 +\tau^2+p^2+3p- 3\tau+3\tau p
    +1\right).
$$
The integration over $\mu$ is over the path $x\mapsto x\eta+\epsilon$,
$x\in\R$ for any generic real $\epsilon$.
\end{thm}

This theorem is proved in \ref{ss-modular}
and \ref{new-modular}.

\section{The integral transformation}\label{s-intt}
Throughout this section, we assume, for definiteness, that
$\mathrm{Im}\,{}\eta<0$, and that $\eta$ is sufficiently small. The results hold for a more general
range of parameters by analytic continuation.

\subsection{A space of functions on which the integral transform is 
  defined} The qKZB heat equation is based on an integral
transformation.  In this section we give a space on which this
integral transformation is defined and invertible.

Let us fix our parameters $\Lambda_i,\eta,\tau$. Then the Shapovalov
form has poles at the points $-\eta\sigma_j$ where
\begin{equation}\label{eq-sigj}
  \sigma_{j}=2(\sum_{k\leq j}\Lambda_k-1),\qquad j=0,\dots,n,
\end{equation}
as well as at the translates of this points by the lattice of periods.

At these points, the hypergeometric integrals obey the following
``resonance relations'' \cite{FV3}.

\medskip

\begin{proposition}\label{p-resrel}
  Let $r,s$, $1\leq a,b\leq n$ be integers and let
  $\sigma_j=2\sum_{k\leq j}\Lambda_k-2$. Then:
\begin{enumerate} 
\item[(i)] If $a<n$, then
\[
u_{a+1,b}(\zz,\eta\sigma_a+r+s\tau,\mu,\tau,p,\eta) = e^{2\pi
  is(z_{a+1}-z_a+\eta\Lambda_{a+1}+\eta\Lambda_a)}
u_{a,b}(\zz,\eta\sigma_a+r+s\tau,\mu,\tau,p,\eta).
\]
\item[(ii)] $u_{1,b}(\zz,-2\eta+r+s\tau,\mu,\tau,p,\eta) = e^{2\pi
    is(z_1-z_n+\eta\Lambda_1+\eta\Lambda_n-p)}
  u_{n,b}(\zz,2\eta+r+s\tau,\mu,\tau,p,\eta).  $
\item[(iii)]If $b<n$, then
\[
u_{a,b+1}(\zz,\lambda,\!-\!\eta\sigma_b+r+sp,\tau,p,\eta) = e^{2\pi
  is(z_{b+1}\!-\!z_b\!-\!\eta\Lambda_{b+1}-\eta\Lambda_b)}
u_{a,b}(\zz,\lambda,\!-\!\eta\sigma_b+r+sp,\tau,p,\eta).
\]
\item[(iv)] $u_{a,1}(\zz,\lambda,2\eta+r+sp,\tau,p,\eta) = e^{2\pi
    is(z_1-z_n-\eta\Lambda_1-\eta\Lambda_n+\tau)}
  u_{a,n}(\zz,\lambda,-2\eta+r+sp,\tau,p,\eta).  $
\end{enumerate}
\end{proposition}

\noindent{\it Proof:} (i) Using the functional relation
$\theta(x+r+s\tau,\tau)=(-1)^{r+s}\exp(-\pi
is(2x+s\tau))\theta(x,\tau)$, we see that
\[
\omega_{a}(t,\zz,\eta\sigma_a+r+s\tau,\tau,\eta)=(-1)^{r+s}e^ {-2\pi i
  s(t-z_a+\eta\Lambda_a)-\pi i s^2\tau}\prod_{j=1}^a\frac
{\theta(t-z_j+\eta\Lambda_j)} {\theta(t-z_j-\eta\Lambda_j)}\,,
\]
and
\[
\omega_{a+1}(t,\zz,\eta\sigma_a+r+s\tau,\tau,\eta)=(-1)^{r+s}e^ {-2\pi
  i s(t-z_{a+1}-\eta\Lambda_{a+1})-\pi i s^2\tau}\prod_{j=1}^a\frac
{\theta(t-z_j+\eta\Lambda_j)} {\theta(t-z_j-\eta\Lambda_j)}\,.
\]
Thus we have equality at the level of integrands.

(ii) Similarly, we have
\[
\omega_{1}(t,\zz,-2\eta+r+s\tau,\tau,\eta)=(-1)^{r+s}e^ {-2\pi i
  s(t-z_1-\eta\Lambda_1)-\pi i s^2\tau},
\]
and using the relation $\sum\Lambda_i=2$ we obtain
\[
\omega_{n}(t,\zz,2\eta+r+s\tau,\tau,\eta)=(-1)^{r+s}e^ {-2\pi i
  s(t-z_{n}+\eta\Lambda_{n})-\pi i s^2\tau}\prod_{j=1}^n\frac
{\theta(t-z_j+\eta\Lambda_j)} {\theta(t-z_j-\eta\Lambda_j)}\,.
\]
The last product of ratios of theta functions may be absorbed using
the functional relation \Ref{eq-zwz} for $\Omega$. This gives
\begin{eqnarray*}
  \lefteqn{ u_{n,b}(\zz,2\eta+r+s\tau,\mu,\tau,p,\eta)
    =e^{-\frac{i\pi}{2\eta}(2\eta+r+s\tau)\mu} } \\ & &\times\int
  e^{-4\pi i\eta} \prod_j\Omega_{\eta\Lambda_j}
  (t-z_j+p,\tau,p)(-1)^{r+s} e^{-2\pi is(t-z_n+\eta\Lambda_n)-\pi i
    s^2\tau} \omega^\vee_b(t,\zz,\mu,p,\eta)dt.
\end{eqnarray*}
The result is then obtained by shifting the integration variable $t$
by $-p$ and using the relation
\[
\omega^\vee_b(t-p,\zz,\mu,p,\eta)=e^{2\pi i(\mu+2\eta)}
\omega^\vee_b(t,\zz,\mu,p,\eta).
\]
The remaining claims (iii), (iv) are proved in a similar way.
\hfill$\square$.

\medskip

\noindent{\bf Definition.} 
Let $E^0(\zz,\tau,c;\eta)$ be the space of holomorphic functions
$\varphi:\C\to V_{\LLambda}[0]$, $\varphi(\lambda)=\sum
\varphi_a(\lambda)\epsilon_a$ such that
\begin{enumerate}
\item[(i)] If $a<n$, then $\varphi_{a+1}(\eta\sigma_a+r+s\tau) =
  e^{2\pi is(z_{a+1}-z_a+\eta\Lambda_{a+1}+\eta\Lambda_a)}
  \varphi_{a}(\eta\sigma_a+r+s\tau).  $
\item[(ii)] $\varphi_{1}(-2\eta+r+s\tau) = e^{2\pi
    is(z_1-z_n+\eta\Lambda_1+\eta\Lambda_n+c)}
  \varphi_{n}(2\eta+r+s\tau). $
\end{enumerate}
Let $E(\zz,\tau,c;\eta)$ be the space of functions $\varphi\in
E^0(\zz,\tau,c;\eta)$ such that
\begin{enumerate}
\item[(iii)] There exist constants $C_1,C_2>0$ (depending on
  $\varphi$) such that
\[
|\varphi_a(\lambda)|\leq C_1\exp
\left(\pi\frac{(\mathrm{Im}\,\lambda)^2}
  {\mathrm{Im}\,{}\tau}+C_2|\lambda| \right),
\]
for all $a=1,\dots,n$.
\end{enumerate}

Examples of functions in $E(\zz,\tau,c;\eta)$ can be constructed using
the universal hypergeometric function:

\medskip

\begin{proposition}\label{p-uinE0} If $\mathrm{Im}\,{}c<0$ then
  the function
\[
\lambda\mapsto \sum_a u_{a,b}(\zz,\lambda,\mu,\tau,-c,\eta)\epsilon_a,
\]
belongs to $E(\zz,\tau,c;\eta)$ for all values of the remaining
parameters.  If $\mathrm{Im}\,{}c>0$ then the function
\[
\mu\mapsto \sum_b u_{a,b}(\zz,\lambda,\mu,c,\tau,-\eta)\epsilon_b,
\]
belongs to $E(\zz,\tau,c;\eta)$ for all values of the remaining
parameters.
\end{proposition}

\medskip

\noindent{\it Proof:} 
Prop.\ \ref{p-resrel} implies that these functions belong to
$E^0(\zz,\tau,c;\eta)$. The bound (iii) follows from Lemma \ref{l-est}
below.  \hfill$\square$

\medskip

\begin{lemma}\label{l-925}
  $ \varphi(\lambda) \in E^0(\zz,\tau, c;\eta)$ if and only if $e^{\pi
    i \lambda^2/4\eta} \varphi(\lambda) \in E^0(\zz,\tau,
  c-\tau;\eta)$ .
\end{lemma}

\medskip

\noindent{\it Proof:}
It is clear that $e^{\frac{i\pi\lambda^2}{4\eta}}\phi(\lambda)$ obeys
(i) in the definition of $E^0(\zz,\tau, c-\tau;\eta)$.  Property (ii)
follows from the identity
\[
e^{\frac {\pi i}{4\eta}(-2\eta+r+s\tau)^2} = e^{-2\pi i s\tau}
e^{\frac {\pi i}{4\eta}( 2\eta+r+s\tau)^2},
\]
for $r,s\in\Z$.  \hfill$\square$

\medskip

For the construction of the integral transform, we will need the
following variant of Prop.\ \ref{p-uinE0}.

\medskip
 
\begin{corollary}\label{c-uinE}
  The function
\[
\lambda\mapsto e^{-\pi i\frac{\lambda^2+\mu^2}{4\eta}} \sum_a
u_{a,b}(\zz,\lambda,\mu,\tau,p,\eta)\epsilon_a,
\]
belongs to $E^0(\zz,\tau,\tau-p;\eta)$ for all $\mu$.  The function
\[
\mu\mapsto e^{\pi i\frac{\lambda^2+\mu^2}{4\eta}} \sum_b
u_{a,b}(\zz,\lambda,\mu,\tau,p,-\eta)\epsilon_b,
\]
belongs to $E^0(\zz,p,\tau-p;\eta)$ for all $\lambda$.
\end{corollary}
\medskip

\noindent{\it Proof:} An immediate consequence of Prop.\ \ref{p-uinE0}
and Lemma \ref{l-925}.  \hfill$\square$

\subsection{The definition of the integral transform}

We define two integral transforms $\F_{\zz,p,\tau}$ and $\tilde
\F_{\zz,\tau,p}$, then show that they are inverse to each other. The
first is $\varphi\mapsto \F_{\zz,p,\tau}(\varphi)=\hat\varphi$ with
\begin{equation}\label{eq-F1}
  \hat\varphi_c(\mu)= \frac {e^{\pi i\frac{\mu^2}{4\eta}}}
  {16\pi^2\eta}\int_{\R+i\sigma} \sum_{b=1}^n \varphi_b(-\nu)
  Q_{b}(\nu,\tau,\eta) u_{b,c}(\zz,\nu,\mu,\tau,p,-\eta) e^{\pi
    i\frac{\nu^2}{4\eta}}d\nu.
\end{equation}
The second is $\psi\mapsto \tilde\F_{\zz,\tau,p}(\psi)=\check\psi$
with
\begin{equation}\label{eq-F2}
  \check\psi_a(\lambda) =e^{-\pi i\frac{\lambda^2}{4\eta}}
  \int_{\eta\R+\tilde\sigma} \sum_{c=1}^n
  u_{a,c}(\zz,\lambda,\mu,\tau,p,\eta) Q_{c}(\mu,p,\eta) \psi_c(-\mu)
  e^{-\pi i\frac{\mu^2}{4\eta}} d\mu.
\end{equation}
These transformations depend on the choice of a real parameter
$\sigma$ or $\tilde\sigma$. We say that $\sigma$ is an {\em admissible
  shift} for $\F$ if it does not lie in any interval
$[\mathrm{Im}(2\eta+ sp), \mathrm{Im}(-2\eta+ sp)]$, $s\in\Z$.

We say that $\tilde\sigma$ is an {\em admissible shift} for $\tilde\F$
if it is a generic real number.

\medskip

\begin{thm}\label{thm-Fourier}
  Suppose that $\mathrm{Im}\,{}\eta<0$.  If $\varphi\in
  E(\zz,\tau,\tau-p;\eta)$ then the integral \Ref{eq-F1} is absolutely
  convergent, independent of the choice of admissible shift and
  defines a linear map $\F_{\zz,p,\tau}:E(\zz,\tau,\tau-p;\eta)\to
  E(\zz,p,\tau-p;\eta)$.  If $\psi\in E(\zz,p,\tau-p;\eta)$ then the
  integral \Ref{eq-F2} is absolutely convergent, independent of the
  choice of admissible shift and defines a linear map
  $\tilde\F_{\zz,\tau,p}:E(\zz,p,\tau-p;\eta)\to
  E(\zz,\tau,\tau-p;\eta)$.  Moreover
\[
\F_{\zz,p,\tau}\circ \tilde\F_{\zz,\tau,p}
=\mathrm{Id}_{E(\zz,p,\tau-p;\eta)}, \qquad \tilde\F_{\zz,\tau,p}\circ
\F_{\zz,p,\tau} =\mathrm{Id}_{E(\zz,\tau,\tau-p;\eta)},
\]
\end{thm}

\medskip

The proof of the first part of the theorem is contained in the next
subsection. The proof of the inversion formula is deferred to
\ref{ss-inv}.

\medskip \newcommand{\zzo}{\zz\,{}^0}

This theorem implies a completeness result for the qKZB equations.
For generic $\zzo\in\C^n$, we may consider the qKZB equations on the
set $\zzo+(p\Z)^n$.  Any solution is uniquely determined by its
initial condition at $\zzo$.  A class of solutions is given by taking
linear combinations of components of the fundamental hypergeometric
solution: Let us say that a solution $v$ is {\em of hypergeometric
  type} if it is of the form:
\[
v_a(z,\lambda) = \int_\gamma \sum_b \bar
u_{ab}(z,\lambda,\mu,\tau,p,\eta) F_b(\mu) d\mu,
\]
for some functions $F_b(\mu)$ and some cycle $\gamma$.

\medskip

\begin{corollary}\label{completeness}  Any solution $v(z,\lambda)$ with
  initial condition in $e^{\frac{\pi i\lambda^2}{4\eta}}
  E(z,\tau,\tau-p;\eta)$ is of hypergeometric type.

  More precisely, let $\zzo\in\C^n$ be generic.  Suppose $\varphi\in
  E(\zz\,{}^0,\tau,\tau-p;\eta)$ and let $\hat\varphi(\mu) =
  \F_{\zzo,p,\tau} (\varphi)(\mu)$.  Then, for all
  $\zz\in\zz\,{}^0+(p\Z)^n$, $\psi(\zz,\mu)=\prod_{i=1}^n
  D_i(-\mu)^{-(z_i-z_i^0)/p}\hat\varphi(\mu)$, viewed as a function of
  $\mu$, belongs to $E(\zz,p,\tau-p;\eta)$ and the function
  $v(\zz,\lambda)=\sum v_a(\zz,\lambda)\epsilon_a$, with
\[
v_a(\zz,\lambda)=\int_{\eta\R+\tilde\sigma}\sum_b
u_{a,b}(\zz,\lambda,\mu,\tau,p,\eta) Q_b(\mu,p,\eta) e^{-\pi
  i\frac{\mu^2}{4\eta}} \psi_b(\zz,-\mu)d\mu,
\]
is, for any generic $\tilde\sigma\in\R$, a solution of hypergeometric
type of the qKZB equations with modulus $\tau$, step $p$ and initial
condition $v(\zz\,{}^0,\lambda)=e^{\pi
  i\frac{\lambda^2}{4\eta}}\varphi(\lambda)$.  Moreover, for all
$\zz\in\zz\,^0+(p\Z)^n$, $v(\zz,\lambda)$ is independent of the choice
of $\tilde\sigma$, and $e^{-i\pi \lambda^2/4\eta} v(\zz,\lambda)$
belongs to $E(\zz,\tau,\tau-p;\eta)$.
\end{corollary}

\medskip

\noindent{\it Proof:} Let $d_{j,a}(\mu)$ denote the eigenvalues
of $D_j(\mu)$:
\[
D_j(\mu)\epsilon_a= d_{j,a}(\mu)\epsilon_a.
\]
They grow at most exponentially as $\mu\to\infty$.  The coordinates of
$\psi$ are then
\[
\psi_a(\zz,\mu)=\prod_{j=1}^n d_{j,a}(-\mu)^{-(z_i-z_i^0)/p}
\hat\varphi_a(\mu).
\]
By the theorem, $\hat\varphi\in E(\zz\,{}^0,p,\tau-p;\eta)$.  Then
$\psi$ clearly obeys the bound (iii) of the definition of
$E(\zz,p,\tau-p;\eta)$. The resonance relations (i),(ii) are checked
by inserting the definitions. For example, for $a<n$ we have (using
the fact that $(\zz-\zzo)/p$ has integral coordinates)
\begin{eqnarray*}
  \psi_{a+1}(\zz,\eta\sigma_a+r+sp) &=&\psi_{a}(\zz,\eta\sigma_a+r+sp)
  e^{2\pi is(z^0_{a+1}-z^0_a+\eta \Lambda_{a+1}+\eta\Lambda_a)} \\ &
  &\times\prod_j \left\{ \frac {d_{j,a+1}(-\eta\sigma_a-r-sp)}
    {d_{j,a}(-\eta\sigma_a-r-sp)} \right\} ^{-\frac{z_j-z_j^0}p}\,.
\end{eqnarray*}
On the other hand, the expression in curly brackets is equal to
$e^{2\pi isp}$ if $j=a$, to $e^{-2\pi isp}$ if $j=a+1$, and to $1$
otherwise.  It then follows by Theorem \ref{thm-Fourier} that $
\check\psi(\zz,\lambda)=e^{-i\pi\lambda^2/4\eta}v(\zz,\lambda)$ is in
$E(\zz,\tau,\tau-p;\eta)$ and that the initial condition
$\check\psi(\zz\,^0) =\varphi(\lambda)$ is satisfied. $v$ is a
solution of hypergeometric type since the $\zz,\lambda$-dependent part
of the kernel of integration is $1\otimes \prod_iD_i(\mu)^{-z_i/p}u$
which is equal to the fundamental solution $\bar u$ up to
$\zz,\lambda$-independent factors, cf.\ \Ref{e-d99}.  \hfill$\square$
\medskip

\noindent{\bf Remarks.}
\begin{enumerate}
\item The result may be expressed in the following terms: our
  generalized Fourier transform maps the qKZB equations to the
  difference equations
\[\psi(\zz+p\delta_i,\mu)=D_i(-\mu)^{-1}\psi(\zz,\mu),\] which are
easily solved, since $D_i$ is a diagonal multiplication operator.  So
the solution given in the corollary is the Fourier transform of the
solution of this simple system of equations with initial condition
given by the inverse Fourier transform of the given initial condition.
\item An initial condition at $\zzo$ uniquely determines a solution
  only on the set $\Gamma=\zzo+(p\Z)^n$, since solutions defined for
  all $\zz$ can always be changed by multiplying by $p$-periodic
  functions without changing the initial condition.  However the
  formula in Corollary \ref{completeness} gives a solution for all
  $\zz$.  On the other hand it is only for $\zz\in \zzo+(p\Z)^n$ that
  we know that the integral is independent of the choice of
  $\tilde\sigma$. For general $\zz$ there is no cancellation of
  residues at pairs of poles that would allow us to move the
  integration contours, and one should expect that the solution
  depends on $\tilde\sigma$. In fact, if we move $\tilde\sigma$ just a
  little bit, then the integration contour crosses infinitely many
  poles (for generic $\eta$), so one would expect that the solution
  depends discontinuously on $\tilde\sigma$. It would be interesting to
  understand this dependence in more detail.
\end{enumerate}

\subsection{Proof of Theorem \ref{thm-Fourier}}

We prove here that the integral transformations are well-defined on
the considered spaces.

We start by estimating the integrands.

\medskip

\begin{lemma}\label{l-est}
\[
|e^{\pi i\frac{\lambda\mu}{2\eta}}
u_{a,b}(\zz,\lambda,\mu,\tau,p,\eta)| \leq C_1 \exp \left(\pi
  \frac{(\mathrm{Im}\,\lambda)^2}
  {\mathrm{Im}\,{}\tau}+C_2|\mathrm{Im}\,{}\lambda| +\pi
  \frac{(\mathrm{Im}\,\mu)^2}{\mathrm{Im}\,{}p}+
  C_2|\mathrm{Im}\,{}\mu|\right),
\]
for some $C_1,C_2>0$ depending on $\zz,\tau,p,\eta,\LLambda$.

For every $\epsilon>0,\zz,\tau,p,\eta,\LLambda$ there exist constants
$C_3,C_4>0$ such that if the distance between $\lambda$ and the
singularities of $Q_a$ is at least $\epsilon$, then
\[
|Q_a(\lambda,\tau,\eta)| \leq C_3 \exp \left(-2\pi
  \frac{(\mathrm{Im}\,\lambda)^2} {\mathrm{Im}\,{}\tau}+
  C_4|\mathrm{Im}\,{}\lambda|\right).
\]
\end{lemma}

\medskip

\noindent{\it Proof:}
The first bound is obtained by applying the estimates of Lemma
\ref{l-C1} to the integral \Ref{eq-467}.  The second follows using the
lower bound in Lemma \ref{l-C1}.  \hfill$\square$

\medskip

In order to show that the integral is independent of the choice of
admissible shift and to bound the integral transform, we will need to
shift the integration contour.  The next result shows that this is
possible since poles occur in pairs with opposite residues.

\medskip

\begin{lemma}\label{l-res}
\begin{enumerate}
\item[(i)] Suppose $\psi=\sum\psi_a\epsilon_a\in
  E(\zz,p,\tau-p;\eta)$. Then for all $r,s\in\Z$, $c=1,\dots,n$ and
  $a=1,\dots,n-1$,
\begin{eqnarray*}
  \lefteqn{ \mathrm{res}_{\mu=-\eta\sigma_a+r+sp} e^{-\pi
      i\frac{\lambda^2+\mu^2}{4\eta}}
    u_{c,a}(\zz,\lambda,\mu,\tau,p,\eta) Q_a(\mu,p,\eta) \psi_a(-\mu)
    } \\ &=& -\mathrm{res}_{\mu=-\eta\sigma_a+r+sp} e^{-\pi
    i\frac{\lambda^2+\mu^2}{4\eta}} u_{c,a+1}(\zz,
  \lambda,\mu,\tau,p,\eta) Q_{a+1}(\mu,p,\eta) \psi_{a+1}(-\mu), \\ 
  \lefteqn{ \mathrm{res}_{\mu=-2\eta+r+sp} e^{-\pi
      i\frac{\lambda^2+\mu^2}{4\eta}}
    u_{c,n}(\zz,\lambda,\mu,\tau,p,\eta) Q_n(\mu,p,\eta) \psi_n(-\mu)
    } \\ &=& -\mathrm{res}_{\mu=2\eta+r+sp} e^{-\pi
    i\frac{\lambda^2+\mu^2}{4\eta}} u_{c,1}(\zz,
  \lambda,\mu,\tau,p,\eta) Q_{1}(\mu,p,\eta) \psi_{1}(-\mu).
\end{eqnarray*}
\item[(ii)] Suppose $\varphi=\sum\varphi_a\epsilon_a\in
  E(\zz,\tau,\tau-p;\eta)$. Then for all $r,s\in\Z$, $c=1,\dots,n$ and
  $a=1,\dots,n-1$,
\begin{eqnarray*}
  \lefteqn{ \mathrm{res}_{\lambda=-\eta\sigma_a+r+s\tau} e^{\pi
      i\frac{\lambda^2+\mu^2}{4\eta}}
    u_{a,c}(\zz,\lambda,\mu,\tau,p,-\eta) Q_a(\lambda,\tau,\eta)
    \varphi_a(-\lambda) } \\ &=&
  -\mathrm{res}_{\lambda=-\eta\sigma_a+r+s\tau} e^{\pi
    i\frac{\lambda^2+\mu^2}{4\eta}} u_{a+1,c}(\zz,
  \lambda,\mu,\tau,p,-\eta) Q_{a+1}(\lambda,\tau,\eta)
  \varphi_{a+1}(-\lambda), \\ \lefteqn{
    \mathrm{res}_{\lambda=-2\eta+r+s\tau} e^{\pi
      i\frac{\lambda^2+\mu^2}{4\eta}}
    u_{n,c}(\zz,\lambda,\mu,\tau,p,-\eta) Q_n(\lambda,\tau,\eta)
    \varphi_n(-\lambda) } \\ &=& -\mathrm{res}_{\lambda=2\eta+r+s\tau}
  e^{\pi i\frac{\lambda^2+\mu^2}{4\eta}} u_{1,c}(\zz,
  \lambda,\mu,\tau,p,-\eta) Q_{1}(\lambda,p,\eta)
  \varphi_{1}(-\lambda).
\end{eqnarray*}
\end{enumerate}
\end{lemma}

\medskip

\noindent{\it Proof:} A straightforward calculation.
\hfill$\square$

\medskip

\noindent{\it Proof of the first part of Theorem \ref{thm-Fourier}.}
It follows from Lemma \ref{l-res} (iii) that the integrand in
\Ref{eq-F1} has only simple poles at $\pm2\eta+\Z+\tau\Z$. If $\sigma$
is an admissible shift then the integration contour stays away from
the singularities and we can use the estimates of Lemma \ref{l-est}.
The integrand is then the product of a function growing at most like
$e^{C|\nu|}$ times $e^{\pi i\nu^2/4\eta}$ which for
$\mathrm{Im}\,{}\eta<0$ converges very rapidly to zero in the real
direction.  So the integrand is an $L^1$ function.  Moreover by Lemma
\ref{l-res} (iv) the residue at $2\eta+r+s\tau$ is opposite to the
residue at $-2\eta+r+s\tau$. It follows that the integration contour
can me moved across each of these pairs of singularities without
changing the value of the integral. This shows that the integral is
independent of the choice of admissible shift.

Let us next show that $\hat\varphi$ defined by \Ref{eq-F1} belongs to
$E(\zz,p,\tau-p;\eta)$.  Properties (i) (ii) follow from Corollary
\ref{c-uinE}, and we are left with the proof of (iii).  For this we
shift variables: let $u^0_{a,b}(\zz,\lambda,\mu,\tau,p,\eta)
=e^{i\pi\frac{\lambda\mu}{4\eta}}u_{a,b}(\zz,\lambda,\mu,\tau,p,\eta)$.
Then
\[
\hat\varphi_c(\mu)= \frac1 {16\pi^2\eta}\int_{\R+i\sigma} \sum_{b=1}^n
\varphi_b(-\nu+\mu) Q_{b}(\nu-\mu,\tau,\eta)
u^0_{b,c}(\zz,\nu-\mu,\mu,\tau,p,-\eta) e^{\pi
  i\frac{\nu^2}{4\eta}}d\nu.
\]
The integration contour is moved by this change of variable, but using
the independence on the choice of $\sigma$ we may choose it to lie in
the interval $[0,\mathrm{Im}\,\tau]$ say.

Then Lemma \ref{l-est} yields
\[
|\hat\varphi_c(\mu)|\leq
\frac{e^{\pi\frac{(\mathrm{Im}\,\mu)^2}{\mathrm{Im}\,{}\tau}}}
{16\pi^2|\eta|}\int_{\R+i\sigma} \sum_{b=1}^n Ce^{C'|\mu-\nu|} |e^{\pi
  i\frac{\nu^2}{4\eta}}|d\nu.
\]
The triangle inequality $|\mu-\nu|\leq |\mu|+|\nu|$ and the Gaussian
integral over $\nu$ give the estimate (iii) in the definition of
$E(\zz,p,\tau-p;\eta)$.

The proof that \Ref{eq-F2} defines a function in
$E(\zz,\tau,\tau-p;\eta)$ is analogous with a little difference. Now
the integration contour is parallel to the line $2\eta\R$. So if the
shift $\tilde\sigma$ is generic it will not meet the singular set
$S=\{\pm2\eta\}+\Z+p\Z$ of the integrand. However it will come
arbitrarily close to these singularities.  Still, the integral is
absolutely convergent, since the distance to the singular set
decreases polynomially: dist$(\mu,S) \geq
\mathrm{const}(1+|\mu|)^{-\alpha}$ for some $\alpha>0$.  This implies
that the divergence coming from the poles close to the integration
contour is at most polynomial. This does not spoil the integrability
which is due to the exponential decay of $\exp(-i\pi\mu^2/4\eta)$.
\hfill$\square$

\medskip

\section{Calculations}
\subsection{The heat equation}\label{ss-he}
Here we prove Theorem \ref{t-HE1}.

The statement of the theorem is
\begin{eqnarray}\label{claim1}
  \lefteqn{u_{a,b}(\zz,\lambda,\nu,\tau,p,\eta)\nonumber =
    -\frac{e^{4\pi i\eta}}{2\pi\sqrt{4i\eta}}
    e^{-i\frac{\pi}{4\eta}(\lambda^2+\nu^2)} } \\ & &\times
  \int_{\eta\R+\tilde\sigma}\sum_{c=1}^n
  u_{a,c}(\zz,\lambda,\mu,\tau,p+\tau,\eta)Q_c(\mu,\tau+p,\eta)
  e^{-\frac{i\pi\mu^2}{4\eta}}
  u_{c,b}(\zz,-\mu,\nu,\tau+p,p,\eta)d\mu.
\end{eqnarray}
We proceed to evaluate the right-hand side. The $\mu$-dependent part
may be simplified by using the following identity:
\begin{lemma}\label{l-89}
\begin{eqnarray*}
  \lefteqn{\sum_{c=1}^n\omega^\vee_c(t,\zz,\mu,\tau+p,\eta)
    Q_c(\mu,\tau+p,\eta)\omega_c(s,\zz,-\mu,\tau+p,\eta) =} \\ &
  &\frac{\theta(\mu-2\eta+t-s,\tau+p)\theta'(0,\tau+p)}
  {\theta(\mu-2\eta,\tau+p)\theta(t-s,\tau+p)} \prod_{j=1}^{n} \frac
  {\theta(t-z_j+\eta\Lambda_j,\tau+p)}
  {\theta(t-z_j-\eta\Lambda_j,\tau+p)} \\ &-&
  \frac{\theta(\mu+2\eta+t-s,\tau+p)\theta'(0,\tau+p)}
  {\theta(\mu+2\eta,\tau+p)\theta(t-s,\tau+p)} \prod_{j=1}^{n} \frac
  {\theta(s-z_j+\eta\Lambda_j,\tau+p)}
  {\theta(s-z_j-\eta\Lambda_j,\tau+p)}\,.
\end{eqnarray*}
\end{lemma}
\noindent{\it Proof:}
Consider each term in the sum on the left as a function of $\mu$: it
is periodic with period one and as $\mu$ is replaced by $\mu+\tau+p$,
it is multiplied by $\exp(-2\pi i(t-s))$. The poles of the term
labeled by $c$ are at $\mu=\mu_c=2\eta-2\eta\sum_{j<c} \Lambda_j$ and
at $\mu=\mu_c-2\eta\Lambda_c=\mu_{c+1}$ modulo $\Z+(\tau+p)\Z$. The
most general form of a meromorphic function of $\mu$ with these
properties is
\[
A_c\,\frac{\theta(\mu-\mu_c+t-s,\tau+p)\theta'(0,\tau+p)}
{\theta(\mu-\mu_c,\tau+p)\theta(t-s,\tau+p)} +B_c\,
\frac{\theta(\mu-\mu_{c+1}+t-s,\tau+p)\theta'(0,\tau+p)}
{\theta(\mu-\mu_{c+1},\tau+p)\theta(t-s,\tau+p)} \, .
\]
The coefficients are determined by comparing the residues at the
poles:
\begin{eqnarray*}
  A_c&=&\omega^\vee_c(t,\zz,\mu_c,\tau+p,\eta)\omega_c(s,\zz,-\mu_c,\tau+p,\eta)
  \\ &=& \prod_{j=c}^{n} \frac {\theta(t-z_j+\eta\Lambda_j,\tau+p)}
  {\theta(t-z_j-\eta\Lambda_j,\tau+p)} \prod_{j=1}^{c-1} \frac
  {\theta(s-z_j+\eta\Lambda_j,\tau+p)}
  {\theta(s-z_j-\eta\Lambda_j,\tau+p)}\,.
\end{eqnarray*}
Similarly one determines $B_c$ which turns out to be equal to
$-A_{c+1}$.  It follows that in the sum over $c$ only two terms are
not canceled and we obtain our claim.  \hfill$\square$

The products of ratios of theta functions in the above identity may be
absorbed into the phase functions $\Omega_{\eta\Lambda_j}$ by means of
their functional relation \Ref{eq-FROmega}. The right-hand side of
\Ref{claim1} is then
\begin{eqnarray}\label{eq-88}
&&{ \frac{-e^{-\frac{i\pi}{2\eta}\lambda\nu}
      }{2\pi\sqrt{4i\eta}} \int
    e^{-\frac{i\pi}{4\eta}(\lambda-\nu+\mu)^2}
    \omega_a(t,\zz,\lambda,\tau,\eta) \omega^\vee_b(s,\zz,\nu,p,\eta)}
  \\ & & \biggl[ \prod_{j=1}^n
  \left(\Omega_{\eta\Lambda_j}(t-z_j+\tau,\tau,\tau+p)
    \Omega_{\eta\Lambda_j}(s-z_j,\tau+p,p)\right) \textstyle
  {\frac{\theta(\mu-2\eta+t-s,\tau+p)\theta'(0,\tau+p)}
    {\theta(\mu-2\eta,\tau+p)\theta(t-s,\tau+p)} } \nonumber \\ & &
  -\prod_{j=1}^n \left(\Omega_{\eta\Lambda_j}(t-z_j,\tau,\tau+p)
    \Omega_{\eta\Lambda_j}(s-z_j+p,\tau+p,p)\right) \textstyle{
    \frac{\theta(\mu+2\eta+t-s,\tau+p)\theta'(0,\tau+p)}
    {\theta(\mu+2\eta,\tau+p)\theta(t-s,\tau+p)} }
  \biggr]dt\,ds\,d\mu.  \nonumber
\end{eqnarray}
The next thing to notice is that if we change variables in the first
of these two terms by replacing $t$ by $t-\tau$, $s$ by $s+p$ and
$\mu$ by $\mu+4\eta$ we obtain exactly the second term up to a sign!
Indeed, the shift of $t$ in $\omega_a$ produces a factor $e^{2\pi
  i(\lambda+2\eta)}$, the shift of $s$ in $\omega^\vee_b$ produces a
factor $e^{-2\pi i(\nu+2\eta)}$ and we get an additional $e^{2\pi
  i\mu}$ from shifting the argument $t-s$ in the ratio of theta
function in the square bracket. These factors are canceled by the
shift of $\mu$ in the exponential function.

To compute this integral we therefore have to carefully consider the
deformation of integration contours involved in the change of
variables.  As we shall see, this deformation produces a residue at
the pole $t=s$.

For these considerations we assume that Im$(\eta\Lambda_j)>N$,
$j=1,\dots,n$, for some $N>0$ large compared to $\tau$, $p$, and that
the points $z_j$ are on the real axis. The general case can then be
obtained by analytic continuation. In this range of parameters, the
integration cycles for the $t$ and $s$ integration is the interval
$[0,1]$. The integrand, viewed as a function of $t$ or $s$, is then
regular in the strip Im$(t)$,Im$(s)\in(-N,N)$.  The first term in
\Ref{eq-88}, however, has additional poles at
$t=s+\alpha+\beta(\tau+p)$, ($\alpha,\beta\in\Z$).  To deal with these
poles we move slightly the $s$-integration cycle into the upper half
plane. After the change of variable, in the first term in \Ref{eq-88},
$\bar t=t+\tau$, $\bar s= s-p$, the new variables are integrated over
$\bar t\in\tau+[0,1]$, $\bar s\in -p+i\epsilon+[0,1]$ for some small
$\epsilon>0$. The first term becomes then equal to the minus the
second term in \Ref{eq-88} after deforming the integration cycles to
the original position, but during this deformation we encounter the
pole at $t=s$. Therefore we obtain a residue
\begin{eqnarray*}
  \Ref{eq-88}&=& \frac{i}{\sqrt{4i\eta}}
  e^{-\frac{i\pi}{2\eta}\lambda\nu} \int
  e^{-\frac{i\pi}{4\eta}(\lambda-\nu+\mu)^2}
  \omega_a(s,\zz,\lambda,\tau,\eta) \omega^\vee_b(s,\zz,\nu,p,\eta) \\ 
  & &\times\prod_{j=1}^n \Omega_{\eta\Lambda_j}(s-z_j,\tau,\tau+p)
  \Omega_{\eta\Lambda_j}(s-z_j+p,\tau+p,p) ds\,d\mu.  \nonumber
\end{eqnarray*}
Using the identity \Ref{eq-Om1} and
\[
\frac{i}{\sqrt{4i\eta}}\int_{\eta\R}
e^{-\frac{i\pi}{4\eta}(\lambda-\nu+\mu)^2}d\mu=1, \qquad
\mathrm{Im}\,{}\eta<0,
\]
we see that this expression reduces to
$u(\zz,\lambda,\nu,\tau,p,\eta)$, completing the proof of Theorem
\ref{t-HE1}.

\subsection{First modular relation}
\label{ss-modular}
Here we present the proof of Theorem \ref{t-mod} (i).
The assumption that $\mathrm{Im}(\eta\tau/p)<0$ implies that the
integration over $\mu$ on the path $x\mapsto \eta x+\epsilon$,
$x\in\R$, converges absolutely: the singularities of $Q$ are avoided
for generic $\epsilon$, and at infinity the factor
$\exp(-i\pi\tau\mu^2/4\eta p)$ converges very fast to zero.  To get
the identity as stated in the Theorem, one uses the following simple
properties of $\rho_a$:
\begin{eqnarray*}
  \rho_a \left(\textstyle{ \frac \lambda\tau, \frac \zz \tau,
      -\frac1\tau, \frac\eta \tau}\right)
  &=&\rho_a(\lambda,\zz,\tau,\eta)^{-1}, \\ \rho_b^\vee
  \left(\textstyle{ \frac \nu p, \frac \zz p, -\frac\tau p, \frac\eta
      p}\right) \rho_b \left(\textstyle{ -\frac \nu p, \frac \zz p,
      -\frac\tau p, \frac\eta p}\right) &=&e^{\frac{2\pi i}{p\tau}
    (\nu-2\eta+2\eta\sum_1^{b-1}\Lambda_j)
    (\nu-2\eta+2\eta\sum_1^{b}\Lambda_j) }.
\end{eqnarray*}
Let us proceed with the proof.  We insert in the left-hand side the
hypergeometric integral for $u_{a,c}$ and $u_{c,b}$ and call the
integration variables $t$ and $s$, respectively. It will be convenient
to make the change of variables $s\to s/p$.  To apply Lemma \ref{l-89}
we use the transformation properties of weight functions:
\begin{equation}\label{eq-tpwf}
  \omega_c\left(\textstyle{ \frac sp , \frac \zz p , -\frac \mu p ,
      -\frac 1p , \frac \eta p} \right) = \rho_c(-\mu,\zz,p,\eta)^{-1}
  e^{-\frac{\pi i}p(2s-\sum_jz_j\Lambda_j)(\mu-2\eta)}
  \omega_c(s,\zz,-\mu,p,\eta).
\end{equation}
A similar formula involving $\rho_b^\vee$ gives the transformation
behavior of $\omega_b^\vee$. Using Lemma \ref{l-89} and the functional
relation \Ref{eq-zwz} of $\Omega_a$, we get:
\begin{eqnarray}\label{eq-48}
  \lefteqn{\sum_{c=1}^n\int u_{a,c} (\zz,\lambda,\mu,\tau,p,\eta)
    u_{c,b} \left(\textstyle{ \frac \zz p, -\frac\mu p, \frac \nu p,
        -\frac1p, -\frac\tau p, \frac\eta p }\right)} \nonumber \\ &
  &\times Q_c(\mu,p,\eta)\rho_c(-\mu,\zz,p,\eta) e^{-\frac{\pi
      i\tau\mu^2}{4\eta p}}\,d\mu \\ &= & \frac1p\int
  \omega_a(t,\zz,\lambda,\tau,\eta)
  \omega_b^\vee(s/p,\zz/p,\nu/p,-\tau/p,\eta/p)
  e^{-\frac{i\pi}{2\eta}(\lambda-\nu/ p)\mu}\nonumber \\ &
  &\times\left[ \Phi(\mu,t,s)e^{-4\pi
      i\eta}-\Phi(\mu+4\eta,t-\tau,s-\tau) e^{-\frac{2\pi
        i\tau}p(\mu+2\eta)-\frac{4\pi i\eta}p} \right]
  \,dt\,ds\,e^{-\frac{\pi i\tau}{4\eta p}\mu^2} d\mu,\nonumber
\end{eqnarray}
where
\begin{eqnarray*}
  \Phi(\mu,t,s) &=& \frac{\theta(\mu-2\eta+t-s,p)\theta'(0,p)}
  {\theta(\mu-2\eta,p)\theta(t-s,p)} e^{-\frac{\pi
      i}p(2s-\sum_{j=1}^nz_j\Lambda_j)(\mu-2\eta)} \\ &
  &\times\prod_{j=1}^n \Omega_{\eta\Lambda_j} (t-z_j+\tau,p,\tau)
  \Omega_{\frac{\eta}p\Lambda_j} \left(\textstyle{
      \frac{s-z_j}p,-\frac1p,-\frac\tau p}\right).  \nonumber
\end{eqnarray*}
As in the proof of the heat equation, the two terms in this equation
cancel after formally changing variables $t\mapsto t-\tau$, $s\mapsto
s-\tau$, $\mu\mapsto\mu+4\eta$ in the first term. However, we have to
carefully see what happens to the integration cycles after this change
of variables. We consider the region of parameters where $z_j\in\R$,
$\mathrm{Im}(\eta\Lambda_j)>>0$ and
$\pi/2>>\mathrm{arg}(p)>\mathrm{arg}(\tau)>0$.

Then the original integration contours may be chosen to be the
interval $[0,1]$. For $z\in\C$, let us denote by $\gamma_z$ the path
\[
r\mapsto r z,\qquad r\in[0,1].
\]
After the change of variable $s\to s/p$, the $s$ integration contour
becomes the path $\gamma_p$.

\begin{lemma}\label{residue}
  Suppose that $p,\tau$ are complex numbers in the upper half plane.
  Let $f(t,s)$ be a meromorphic function of two variables such that
  $f(t+1,s)=f(t,s+p)=f(t,s)$ and such that $\alpha(t,s)= (t-s)f(t,s)$
  is regular in a domain containing
\[
\{(t+r\tau,sp+r\tau)\,|\, t,r,s\in[0,1]\}.
\]
Then
\[
\int_{\gamma_1\times\gamma_p} f(t+\tau,s+\tau)dt\wedge ds =
\int_{\gamma_1\times\gamma_p} f(t,s)dt\wedge ds+2\pi
i\int_{\gamma_\tau}\alpha(t,t)dt,
\]
and all integrals are absolutely convergent.
\end{lemma}

\noindent{\it Proof:} 
Since $p$ is not real the integration contours intersect
transversally, so that $f(t,s)$ (which has at most a simple pole on
the diagonal) is absolutely integrable. More generally, with our
assumption on the regularity of $f$, the integral
\[
I_r=\int_{\gamma_1\times\gamma_p} f(t+r\tau,s+r\tau)dt\wedge ds
\]
is absolutely convergent for all $r\in [0,1]$. Let, for $\epsilon>0$,
$D_\epsilon$ be the integration domain obtained from
$\gamma_1\times\gamma_p$ by removing the points where
$|t-s|\leq\epsilon$ and set
\[
I_r(\epsilon)=\int_{D_\epsilon} f(t+r\tau,s+r\tau)dt\wedge ds.
\]
Then $I_1-I_0=\lim_{\epsilon\to 0}\int_0^1\frac d{dr}I_r(\epsilon)dr$.
We have
\begin{eqnarray*}
  \frac d{dr}I_r(\epsilon)&=&-\tau\int_{D_{\epsilon}}
  d(f(t+r\tau,s+r\tau)(dt-ds)) \\ &=& -\tau\int_{\partial
    D_{\epsilon}}f(t+r\tau,s+r\tau)(dt-ds) \\ &=& -\tau\int_{\partial
    D_{\epsilon}}\alpha(t+r\tau,s+r\tau) \frac{dt-ds}{t-s} \\ 
  &=&\tau\int_{|u|=\epsilon}\alpha(u+r\tau,r\tau)\frac{du}{u}
  +O(\epsilon) \\ &=&2\pi i\tau\alpha(r\tau,r\tau),
\end{eqnarray*}
from which the claim follows immediately. In the application of
Stokes' theorem above the other boundary components of $\partial
D_\epsilon$ do not contribute since the function is periodic. The
remaining component is homotopic to a circle with negative
orientation, which explains the change of sign.  \hfill$\square$

Applying this lemma to our situation after shifting $\mu$ by $4\eta$
in the first term yields:
\begin{eqnarray*}
  \lefteqn{\Ref{eq-48}=\frac1p\int_{-\eta\infty}^{\eta\infty}\int_{\gamma_\tau}
    \omega_a(t,\zz,\lambda,\tau,\eta)
    \omega_b^\vee(t/p,\zz/p,\nu/p,-\tau/p,\eta/p)
    e^{-\frac{i\pi}{2\eta}(\lambda-\nu/ p)\mu-\frac {4\pi i\eta}{p}}
    }\\ & &\times e^{-\frac{\pi i}p(2t-\sum_jz_j\Lambda_j)(\mu+2\eta)}
  \prod_{j=1}^n \Omega_{\eta\Lambda_j} (t-z_j,p,\tau)
  \Omega_{\frac{\eta}p\Lambda_j} \left(\textstyle{
      \frac{t-z_j-\tau}p,-\frac1p,-\frac\tau p}\right)
  \,dt\,e^{-\frac{\pi i\tau}{4\eta p}\mu^2} d\mu.
\end{eqnarray*}
Now the Gaussian integral over $\mu$ along $\eta\R$ may be performed
explicitly, and our claim follows from \Ref{eq-Om4} and \Ref{eq-tpwf}.

\subsection{Second modular relation}\label{new-modular}
Here we present the proof of Theorem \ref{t-mod}
(ii) which is parallel to
the proof of Theorem \ref{t-mod} (i).

The assumption that $\mathrm{Im}(\eta p/\tau)<0$ implies that the
integration over $\mu$ over the path $x\mapsto \eta x+\epsilon$,
$x\in\R$, converges absolutely: the singularities of $Q$ are avoided
for generic $\epsilon$, and at infinity the factor
$\exp(-i\pi p\mu^2/4\eta \tau)$ converges very fast to zero.  

  We insert in the left-hand side the
hypergeometric integral for $u_{a,c}$ and $u_{c,b}$ and call the
integration variables $t$ and $s$, respectively. 
It will be convenient
to make the change of variables $t\to t/\tau $.  To apply Lemma \ref{l-89}
we use the transformation properties of weight functions:
\begin{equation}\label{new-eq-tpwf}
  \omega_c^\vee\left(\textstyle{ \frac t\tau , \frac \zz \tau , 
\frac \mu \tau  ,
      -\frac 1\tau , \frac \eta \tau} \right) = 
\rho_c^\vee(\mu,\zz,\tau ,\eta)^{-1}
  e^{\frac{\pi i}\tau (2t-\sum_jz_j\Lambda_j)(\mu + 2\eta)}
  \omega_c^\vee(t,\zz,\mu,\tau,\eta).
\end{equation}
We have
\begin{eqnarray}\label{new-proof}
\lefteqn{    \sum_{c=1}^n \int   
u_{a,c}\left(\textstyle{ \frac 
      \zz\tau, \frac\lambda\tau, \frac\mu\tau, -\frac p\tau, 
-\frac1\tau, \frac\eta\tau }\right)
\, u_{c,b}(\zz,-\mu,\nu,\tau,p,\eta) \nonumber Q_c(\mu,\tau,\eta) \,
\rho^\vee_c (\mu,\zz,\tau,\eta)\,
e^{-{i\pi p\over 4\eta \tau}\mu^2}\,d\mu}  \\ 
&= &
\frac1\tau \int
  \omega_a(t/\tau,\zz/\tau,\lambda/\tau,-p/\tau,\eta/\tau)
  \omega_b^\vee(s,\zz,\nu,p,\eta)
  e^{-\frac{i\pi}{2\eta}(\lambda/\tau-\nu)\mu}\nonumber \\ &
 {} &\times\left[ \Phi(\mu,t,s)e^{-{4\pi i\eta\over \tau}}
   -\Phi(\mu+4\eta,t+p,s+p) 
 e^{-\frac{2\pi i p}\tau(\mu -2\eta)-{4\pi i\eta}} \right]
  \,dt\,ds\,e^{-\frac{\pi i p}{4\eta \tau}\mu^2} d\mu,\nonumber
\end{eqnarray}
where
\begin{eqnarray*}
  \Phi(\mu,t,s) &=& \frac{\theta(\mu-2\eta+t-s,\tau)\theta'(0,\tau)}
  {\theta(\mu-2\eta,\tau)\theta(t-s,\tau)} e^{\frac{\pi
      i}\tau (2t-\sum_{j=1}^nz_j\Lambda_j)(\mu+2\eta)} \\ &
  &\times\prod_{j=1}^n \Omega_{{\eta\over \tau}\Lambda_j} 
\left({t-z_j-p\over \tau},-{p\over\tau},-{1\over\tau}\right)
  \Omega_{\eta\Lambda_j} (
s-z_j,\tau,p).  \nonumber
\end{eqnarray*}
The two terms in the right hand side of this equation
cancel after formally changing variables $t\mapsto t-p$, $s\mapsto
s-p$, $\mu\mapsto\mu-4\eta$ in the second term. However, we have to
carefully see what happens to the $(t,s)$ integration
cycle after this change of variables. 

We consider the region of parameters where $z_j\in\R$, 
$\mathrm{Im}(\eta\Lambda_j)>>0$ and
$\pi/2>>\mathrm{arg}(\tau)>\mathrm{arg}(p)>0$.
In this case the $(t,s)$ integration cycle is the product of paths
$\gamma_\tau\times\gamma_1$.

Applying to this situation Lemma \ref{residue} after 
decreasing  $\mu$ by $4\eta$ we see that the left hand side of the new
modular relation is equal to
$$
\frac{2\pi i}\tau\int_{-\eta\infty}^{\eta\infty}\int_{\gamma_p}
  \omega_a(t/\tau,\zz/\tau,\lambda/\tau,-p/\tau,\eta/\tau)
  \omega_b^\vee(t,\zz,\nu,p,\eta)
  e^{-\frac{i\pi}{2\eta}(\lambda/\tau-\nu)\mu}
e^{-\frac{\pi i p}{4\eta \tau}\mu^2}
$$
$$
\times e^{-{4\pi i\eta\over \tau}}
 e^{\frac{\pi i}\tau (2t-\sum_{j=1}^nz_j\Lambda_j)(\mu+2\eta)} 
\prod_{j=1}^n \Omega_{{\eta\over \tau}\Lambda_j} 
\left({t-z_j-p\over \tau},-{p\over\tau},-{1\over\tau}\right)
  \Omega_{\eta\Lambda_j} (s-z_j,\tau,p). 
$$
Now the Gaussian integral over $\mu$ along $\eta\R$ may be performed
explicitly, and our claim follows from \Ref{eq-Om4} and \Ref{eq-tpwf}.

\subsection{The inversion formulae}\label{ss-inv}

Here we conclude the proof of Theorem \ref{thm-Fourier}.  We have to
prove the identities
\begin{equation}\label{eq-claim4}
  \frac{1}{16\pi^2\eta}\int_{\R+i\sigma}\sum_{c=1}^n u_{c,a}
  (\zz,-\mu,\lambda,p,\tau,-\eta) u_{c,b} (\zz,\mu,\nu,p,\tau,\eta)
  Q_c(-\mu,p,\eta) \,d\mu
  =\frac{\delta_{ab}\delta(\lambda+\nu)}{Q_b(\nu,\tau,\eta)},
\end{equation}
\begin{equation}\label{eq-claim4'}
  \frac{1}{16\pi^2\eta}\int_{\eta\R+\tilde\sigma}\sum_{c=1}^n u_{a,c}
  (\zz,\lambda,\mu,\tau,p,\eta) u_{b,c} (\zz,\nu,-\mu,\tau,p,-\eta)
  Q_c(\mu,p,\eta) \,d\mu
  =\frac{\delta_{ab}\delta(\lambda+\nu)}{Q_b(\nu,\tau,\eta)}.
\end{equation}
We give a proof of the first of these identities. The second is
treated in a similar way.

\begin{figure}
  {\scalebox{.3}{\includegraphics{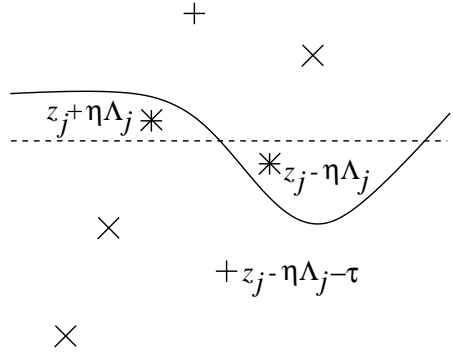}}}
\caption{The integration cycle. The solid line is the contour
  for the integration over $t$ and the dashed line is the contour of
  integration for $s$. They are oriented from left to right.  The
  symbols $\times,+$ indicate the poles of the integrand as a function
  of $t$ and as a function of $s$, respectively. The points
  $z_j\pm\eta\Lambda_j$ are poles for both variables}\label{Fig2}
\end{figure}

As in the proof of the heat equation, the $\mu$-dependent part of the
integrand may be simplified by collecting the terms with the same
poles, and using the functional relation of $\Omega_a$, with the
result:
\begin{eqnarray}\label{eq-90}
  \lefteqn{ \frac1{16\pi^2\eta} \int\sum_{c=1}^n u_{c,a}
    (\zz,-\mu,\lambda,p,\tau,-\eta) u_{c,b} (\zz,\mu,\nu,p,\tau,\eta)
    Q_c(-\mu,p,\eta) \,d\mu}\nonumber\\ &= & \int
  \omega_a^\vee(t,\zz,\lambda,\tau,-\eta)
  \omega_b^\vee(s,\zz,\nu,\tau,\eta)
  e^{-\frac{i\pi}{2\eta}(\lambda+\nu)\mu}\nonumber \\ & &\biggl[ -
  \frac{\theta(\mu+2\eta+s-t,p)\theta'(0,p)}
  {\theta(\mu+2\eta,p)\theta(s-t,p)} \prod_{j=1}^n
  \Omega_{-\eta\Lambda_j} (t-z_j,p,\tau) \Omega_{\eta\Lambda_j}
  (s-z_j,p,\tau) \\ & & +\frac{\theta(\mu-2\eta+s-t,p)\theta'(0,p)}
  {\theta(\mu-2\eta,p)\theta(s-t,p)} \prod_{j=1}^n
  \Omega_{-\eta\Lambda_j} (t\!-\!z_j\!+\!\tau,p,\tau)
  \Omega_{\eta\Lambda_j} (s\!-\!z_j\!+\!\tau,p,\tau)
  \biggr]\,dt\,ds\,d\mu.\nonumber
\end{eqnarray}
As in the proof of the heat equation, the two terms in this equation
cancel after formally changing variables $t\mapsto t-\tau$, $s\mapsto
s-\tau$, $\mu\mapsto\mu+4\eta$ in the second term.  However, we have
to carefully see what happens to the integration cycles after this
change of variables. We consider a region of parameters where
$z_j\in\R$ and $\eta,\tau$ are small compared to 1, Im$(p)$ and to the
distance between the $z_j$.  We study the singularities of the
integrand and the form of the integration cycles in the vicinity of
$z_j$.  The poles of the integrand of $u(\zz,\mu,\nu,p,\tau,\eta)$ are
at $s=z_j-\eta\Lambda_j-k\tau-lp$ and at
$s=z_j+\eta\Lambda_j+k\tau+lp$ ($k,l\in\Z_{\geq 0}$).  The cycle for
the integration over $s$ goes below $z_j+\eta\Lambda_j$ and above
$z_j-\eta\Lambda_j$.  The poles of the integrand of
$u(\zz,-\mu,\lambda,p,\tau,-\eta)$ in the vicinity of $z_j$ are at
$t=z_j+\eta\Lambda_j-k\tau-lp$ and at $t=z_j-\eta\Lambda_j+k\tau+lp$
($k,l\in\Z_{\geq 0}$). The $t$ integration cycle goes above
$z_j+\eta\Lambda_j$ and below $z_j-\eta\Lambda_j$. These integration
contours are depicted in Fig.~\ref{Fig2}.

\begin{figure}
  {\scalebox{.3}{\includegraphics{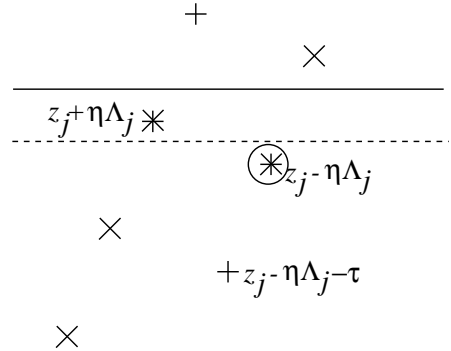}}}
\caption{The integration cycle of Fig.~\ref{Fig2}
  can be deformed to this cycle. The circle around $z_j-\eta\Lambda_j$
  is oriented counterclockwise}
\label{Fig3}
\end{figure}

To treat the two terms on the right-hand side of \Ref{eq-90}
separately, we split the contour for the $t$ integration into two
pieces, as shown in Fig.~\ref{Fig3}. Then the integration cycle does
not meet the additional poles at $t=s$ of both terms.

\begin{figure}
  {\scalebox{.3}{\includegraphics{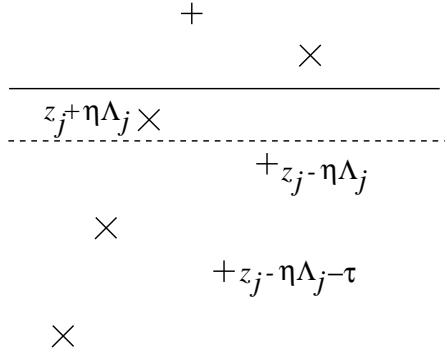}}}
\caption{In the first term of the right-hand side of
  \Ref{eq-90} the integrand is more regular and the cycle may be
  replaced by this one}\label{Fig4}
\end{figure}

Let us first consider the first term in \Ref{eq-90}. As a function of
$t$ it is actually regular at $z_j-\eta\Lambda_j$ so that the circle
around this point does not contribute.  In Fig.~\Ref{Fig4} the
remaining integration cycle is shown, along with the position of the
poles for $t$ and $s$.

\begin{figure}
  {\scalebox{.3}{\includegraphics{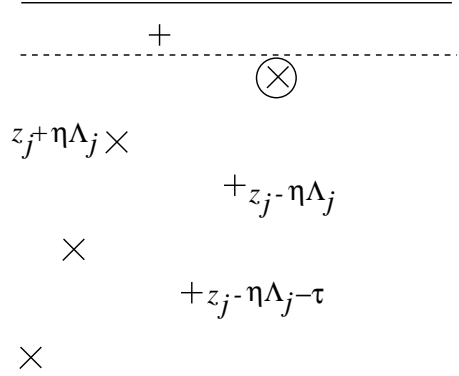}}}
\caption{This is the integration cycle after the shift
  of the variables $t$ and $s$ in the second term.  The poles of the
  result are indicated by the symbols $\times,+$ as in the preceding
  figures}\label{Fig5}
\end{figure}

After changing variables $t\mapsto t-\tau$, $s\mapsto s-\tau$,
$\mu\mapsto\mu+4\eta$ in the second term of \Ref{eq-90} the integrand
of the second term becomes equal to the first one, but the integration
contour are shifted: the $\mu$-integration cycle can be deformed back
to the original one without encountering singularities.  The cycle for
the $t$ and $s$ integration is depicted in Fig.~\ref{Fig5}.

We have to deform the integration cycle to the position shown in
Fig.~\ref{Fig4}. In so doing, we pick the residue at
$t=z_j-\eta\Lambda_j+\tau$ of the residue at $s=t$. At this point, we
may use the identity \Ref{eq-Om5}. The result is
\begin{eqnarray*}
  \Ref{eq-90}&=& \frac{-(2\pi i)^2}{16\pi^2\eta}
  \mathrm{res}_{t=z_j-\eta\Lambda_j+\tau}
  \omega_a^\vee(t,\zz,\lambda,\tau,-\eta)
  \omega_b^\vee(t,\zz,\nu,\tau,\eta) \int
  e^{-\frac{i\pi}{2\eta}(\lambda+\nu)\mu}d\mu \\ &=&
  \mathrm{res}_{t=z_j-\eta\Lambda_j}
  \omega_a^\vee(t,\zz,\lambda,\tau,-\eta)
  \omega_b^\vee(t,\zz,\nu,\tau,\eta) \delta(\lambda+\nu).
\end{eqnarray*}
The residue may be computed at $z_j-\eta\Lambda_j$ since its argument
is $\tau$-periodic.  Our claim follows then from the
\begin{lemma}
\[
\sum_{j=1}^n\mathrm{res}_{t=z_j-\eta\Lambda_j}
\omega_a^\vee(t,\zz,\lambda,\tau,-\eta)
\omega_b^\vee(t,\zz,-\lambda,\tau,\eta)
=\delta_{ab}Q_a(-\lambda,\tau,\eta)^{-1}.\]
\end{lemma}

\noindent{\it Proof:}
This formula can be deduced from the more general result in \cite{TV},
formula (C.4). For the sake of completeness, we include the simple
proof in this special case.  The product of weight functions of which
we compute the residue is an elliptic function of $t$ whose simple
poles in a fundamental domain are in the set $\{z_j\pm\eta\Lambda_j,
j=1,\dots,n\}$.

If $a>b$, then the function is regular at $z_j-\eta\Lambda_j$,
$j=1,\dots,n$, so the residues vanish. If $a<b$, then the function is
regular at $z_j+\eta\Lambda_j$. Therefore the sum of the residues in
the claim is the sum over the residues at {\em all} poles and vanishes
by the residue theorem.

If $a=b$, only the residue at $t=z_a-\eta\Lambda_a$ gives a
non-vanishing contribution, which is easily evaluated, and gives
$Q_a(-\lambda,\tau,\eta)^{-1}$.  \hfill$\square$.

\section{$\SL(3,\Z)$-identities for
  hypergeometric integrals}\label{s-GL}
In this section we recast our results into a form which shows
the analogy with the identities discovered in \cite{FV2} for
the elliptic gamma function. In \cite{FV2} we showed that
the elliptic gamma function is a ``degree 1'' generalized
Jacobi modular function for the group 
$G=\SL(3,\Z)\semidirect\Z^3$. 
This amounts to the identities of Theorem \ref{t-JM} below.
The hypergeometric integral $u$ obeys 
a non-Abelian version of these identities.
We also discuss the geometric interpretation
of these identities: they can be
interpreted as the projective flatness of a discrete connection on
a vector bundle over a suitable $G$-space, see also \cite{FV4}.
As in the case of Gamma functions these properties are more
transparent if we rewrite the identities, which
we formulate in \ref{ss-ide}
 in homogeneous
variables $x$ defined by $\tau=x_1/x_3$, $p=x_2/x_3$. The
$G$-action is then linear. The next step,
as in the case of $\Gamma$, is to extend the
range of parameter to a $G$-space, as the condition of positivity of
imaginary parts is not preserved by the action. This can
be done in a  surprisingly easy way by reflection arguments.
We do it in the simplest case of one tensor
factor in \ref{ss-newrange}. After this is done, we 
rewrite the identities as the projective flatness of
an $\SL(3,\Z)$-connection in \ref{ss-profl}.

\subsection{Identities for the elliptic gamma function}
We first recall the properties of the elliptic gamma
function that are relevant here. More details are included in
Appendix \ref{App-A}.
\begin{thm} \cite{FV1}\label{t-JM}
The elliptic gamma function, defined by the formula
\[\Gamma(z,\tau,\sig)=\prod_{j,k=0}^\infty
\frac
{1-e^{2\pi i((j+1)\tau+(k+1)\sig-z)}}
{1-e^{2\pi i(j\tau+k\sig+z)}}\,,
\]
for $\mathrm{Im}(\tau),\mathrm{Im}(p)>0$,
obeys the following identities
\begin{eqnarray}
  \Gamma(z+1,\tau,\sig)&=&\Gamma(z,\tau,\sig),\label{e-gf2}\\
  \Gamma(z+\sig,\tau,\sig)&=&
  \theta_0(z,\tau)\Gamma(z,\tau,\sig),\label{e-gf3}\\
  \Gamma(z+\tau,\tau,\sig)&=&
  \theta_0(z,\sig)\Gamma(z,\tau,\sig),\label{e-gf4}
\end{eqnarray}
\begin{eqnarray}
\Gamma(z,\tau+1,\sig)&=&\Gamma(z,\tau,\sig+1)=
\Gamma(z,\tau,\sig),\label{e-gf5}
\\
\Gamma(z,\tau+\sig,\sig)
&=&
\frac{\Gamma(z,\tau,\sig)}
{\Gamma(z+\tau,\tau,\sig+\tau)},\label{e-gf6}
\\
\Gamma(
 z/\tau,-1/\tau,\sig/\tau)&=&
e^{i\pi Q(z;\tau,\sig)}           
{\Gamma(
{({z-\tau})/\sig,-\tau/\sig,-1/\sig})}
{\Gamma(z,\tau,\sig)}
,\label{e-gf7}
\end{eqnarray}
for some $Q\in\Q(\tau,\sig)[z]$ given in  Appendix \ref{App-A}.
\end{thm}

In fact $\Gamma$ is also defined for negative imaginary parts
of $\tau$ and $p$ by a reflection procedure,
see \ref{A-ext}, and the above identities continue to hold in this
wider range of parameters.

\subsection{Identities for hypergeometric integrals}\label{ss-ide}

Let $E_\delta(\zz,\tau,c;\eta)$ be the space of functions $\varphi\in
E^0(\zz,\tau,c;\eta)$ such that there exist constants $C_1,C_2>0$
(depending on $\varphi$) such that
\[
|e^{\frac{\pi i\delta\lambda^2}{4\eta}}\varphi_a(\lambda)|\leq C_1\exp
\left(\pi\frac{(\mathrm{Im}\,\lambda)^2}
  {\mathrm{Im}\,{}\tau}+C_2|\lambda| \right),
\]
for all $a=1,\dots,n$.

For $\delta=0$ we have constructed examples of functions in $E_\delta$
in Prop.~\ref{p-uinE0}.

\begin{proposition}
  Let $\Phi_S(\zz,\tau,c,\eta)$, $\Phi_T(\eta)$ be the operators
  acting on $V_\LLambda[0]$-valued functions
  $\varphi(\lambda)=\sum_{a=1}^n\varphi_a(\lambda) \epsilon_a$ as
\begin{eqnarray*}
  (\Phi_S(\zz,\tau,c,\eta)\varphi)_a(\lambda) &=& e^{-\frac{\pi
      ic\lambda^2}{4\eta \tau}}
  \rho_a(\lambda,\zz,\tau,\eta)\varphi_a\left(\frac\lambda\tau\right),
  \\ (\Phi_T(\eta)\varphi)_a(\lambda) &=& e^{-\frac{\pi
      i\lambda^2}{4\eta}} \varphi_a(\lambda), \\ 
  (\Phi_C\varphi)_a(\lambda)&=&\varphi_a(-\lambda).
\end{eqnarray*}
Then these operators restrict to isomorphisms
\begin{eqnarray*}
  \Phi_S(\zz,\tau,c,\eta): E_\delta\left( \frac \zz\tau,  -\frac1\tau,
    \frac c \tau, \frac\eta\tau \right) &\to&
  E_{\frac{\delta+c}\tau}(\zz,\tau,c,\eta), \\ \Phi_T(\eta):
  E_\delta(\zz,\tau,c,\eta) &\to& E_{\delta+1}(\zz,\tau,c+\tau,\eta),
  \\ \Phi_C:E_\delta(-\zz,\tau,-c,-\eta) &\to&
  E_{-\delta}(\zz,\tau,c,\eta).
\end{eqnarray*}
\end{proposition}

Then Theorem \ref{thm-Fourier} about the Fourier transform may be 
reformulated, in a slightly generalized form, as follows:

\medskip

\begin{proposition}\label{invert}
  The operator
\[
(V(\zz,\tau,p,\eta)\varphi)_a(\lambda)= \int
\sum_{b=1}^nu_{ab}(\zz,\lambda,\mu,\tau,p,\eta)Q_b(\mu,p,\eta)\varphi_b(-\mu)d\mu,
\]
is an invertible linear map from $E_\delta(\zz,p,\tau,\eta)$ onto
$E_{-1/\delta}(\zz,\tau,-p,\eta)$.
\end{proposition}

\medskip

The properties of the universal hypergeometric function may then be
expressed as relations for these operators: first of all the analogue
of the identities \Ref{e-gf2}--\Ref{e-gf4} are the qKZB and
mirror qKZB equations \Ref{eq1}. The remaining identities
are given by the following result.

\medskip

\begin{thm}\label{t-relations}
\begin{enumerate}
\item[(i)]
  $E_\delta(\zz,\tau+r,p+s,\eta)=E_\delta(\zz,\tau,p,\eta)$, if
  $r,s\in\Z$.  The identities
\begin{displaymath}
  V(\zz,\tau+1,p,\eta)=V(\zz,\tau,p+1,\eta)=V(\zz,\tau,p,\eta)
\end{displaymath}
hold in $\mathrm{Hom}(E_\delta(\zz,p,-\tau,\eta),
E_{-1/\delta}(\zz,\tau,p,\eta))$.
\item[(ii)] The identity
\begin{displaymath}
  V(\zz,\tau,p,\eta) = c_T \, \Phi_T(\eta) \, V(\zz,\tau,\tau+p,\eta)
  \, \Phi_T(\eta) \, V(\zz,\tau+p,p,\eta) \, \Phi_T(\eta)
  \end{displaymath}
  holds in $\mathrm{Hom}(E_\delta(\zz,p,\tau,\eta),
  E_{-1/\delta}(\zz,\tau,-p,\eta))$ with
\[
c_T=-\frac{e^{4\pi i\eta}} {2\pi\sqrt{4i\eta}}.
\]
\item[(iii)] The identity
\begin{eqnarray*}
  \lefteqn{ \Phi_S(\zz,\tau,-p,\eta) \, V\left(\textstyle{ \frac
        \zz\tau, -\frac1\tau, \frac p\tau, \frac\eta\tau }\right) \,
    \Phi_S\left(\textstyle{ \frac \zz\tau, \frac p\tau, -\frac1\tau,
        \frac\eta\tau }\right) }\\ & &= c_S \, V(\zz,\tau,p,\eta) \,
  \Phi_S(\zz,p,\tau,\eta) \, V\left(\textstyle{ \frac \zz p, -\frac1p,
      -\frac\tau p, \frac\eta p} \right)
\end{eqnarray*}
holds in $\mathrm{Hom}(E_\delta(\zz/p,-\tau/p,-1/p,\eta/p),
E_{\delta'}(\zz,\tau,-p,\eta))$ with $\delta'=p\delta/(1-\tau\delta)$
and
\begin{displaymath}
  c_S=\sqrt{\frac{ip}{4\eta\tau}}\,e^{\frac{\pi i\eta}{3p\tau}
    \left(3\sum_{j<k}\Lambda_j\Lambda_k(z_j-z_k)^2
      +2\left(\sum_{j=1}^n\eta^2\Lambda_j^3 +\tau^2+p^2-3p+3\tau+3\tau
        p +1\right) \right)}.
\end{displaymath}
\end{enumerate}
\end{thm}
 
\medskip

\noindent{\it Proof:}
The first statement (i) is trivial. The second statement (ii) 
is a reformulation of Theorem \ref{t-HE1} and (iii)
is a reformulation  of Theorem \ref{t-mod}. \hfill$\square$

\subsection{The universal hypergeometric function for one tensor factor}\label{ss-onefactor}
For simplicity  we assume
from now on 
  that the tensor product $V_\Lambda$ of $sl_2$ representations consists of
  one factor which is the three-dimensional representation $V_2$. In
  this case the universal hypergeometric function does not depend on
  $z_1$ and is a scalar function,
\begin{equation}\label{eq-u}
  u(\lambda,\mu,\tau,p,\eta) =
  e^{-\frac{i\pi\lambda\mu}{2\eta}} \int_C
  \Omega_{2\eta}(t,\tau,p)
{\theta(\lambda +t,\tau)  \theta(\mu +t,p)  
\over
\theta(t-2\eta,\tau)  \theta(t-2\eta,p)} dt.
\end{equation}
The integration cycle $C$ is as in Sect.\ \ref{s-hyp}:
in this case it can be defined to be the interval $[0,1]$
if $\eta$ has positive imaginary part. For general $\eta$
the integral is defined by analytic continuation, with the
effect of deforming the contour as in Fig.~\ref{Fig1}. 
The resulting function is a meromorphic function
on $\C\times\C\times H_+\times H_+\times\C$, where
$H_+$ denotes the open upper half-plane. This function
is regular as long as the integration cycle is not
pinched between poles of the integrand. As these poles
are at $\pm (2\eta+r\tau+s p+t)$, $r,s\in\Z_{\geq0}, t\in\Z$,
$u$ is regular if $4\eta$ does  not belong to the
set $\Z\tau+\Z p +\Z$.

The following is an interesting alternative formula
for the function $u$,
\begin{eqnarray*}
u(\lambda,\mu,\tau,p,\eta)
&=&-\frac{\Gamma(4\eta,\tau,p)}
{\prod_{j=1}^\infty[(1-e^{2\pi ij\tau})(1-e^{2\pi ijp})]}
e^{-\frac{\pi i}{2\eta}(\lambda+2\eta)(\mu+2\eta)}\\
&&\times\,
\sum_{j,k=0}^\infty
e^{-2\pi i(j\lambda+k\mu)+(2j+1)(2k+1)\eta}
\theta_0(\lambda+2\eta+kp,\tau)
\theta_0(\mu+2\eta+j\tau,p)
\\
&&
\times\,\prod_{l=0}^{k-1}
\frac
{\theta_0(lp+4\eta,\tau)}
{\theta_0((l+1)p,\tau)}
\prod_{l=0}^{j-1}
\frac
{\theta_0(l\tau+4\eta,p)}
{\theta_0((l+1)\tau,p)}\,,
\end{eqnarray*}
where $\theta_0(z,\tau)=\prod_{j=0}^\infty
(1-e^{2\pi i (j\tau+z)})
(1-e^{2\pi i ((j+1)\tau-z)})$, see App.\ \ref{App-A}.
This formula can be proved by moving $C$ to infinity
and picking up the residues at the poles, and is valid
in a certain region of parameter space.
We will not use this formula in this paper.

The Shapovalov form is
\[
Q(\lambda,\tau,\eta) =
\frac{\theta(4\eta,\tau)\theta'(0,\tau)}
{\theta(\lambda-2\eta,\tau)
  \theta(\lambda+2\eta,\tau)}\,.
\]

\subsection{Spaces and operators }
 Consider $\C^3$ and the projectivization of the dual space,
$P(\C^3)^*$. Consider $W \subset (\C-0)^3\times P(\C^3)^*$ where
\[
W=\{((x_1,x_2,x_3),(y_1:y_2:y_3))\,|\, x_1y_1+x_2y_2+x_3y_3=0\,\}\,.
\notag
\]
The natural projection $W\to \C^3-0$ is a projective
 line bundle. The group $SL(3,\Z)$
acts on $W$, $g:(x,y)\mapsto (gx, (g^t)^{-1}y)$ for $g\in
SL(3,\Z)$. 

For $(x;y)\in W$ such that 
$\mathrm{Im}\,{x_1}/{x_3}\neq 0, y_2\neq 0$
introduce the space of functions
$$
E(x;y)=\left\{\begin{array}{rl} 
\{v(\lambda)\,|\,v(x_3\lambda)\in
E_{\delta=-y_1/y_2}(\frac{x_1}{x_3},
\frac{x_2}{x_3},\frac{1}{2x_3})\},&
\mathrm{if}\ \mathrm{Im}\,\frac{x_1}{x_3}> 0,\\
\{v(\lambda)\,|\, v(x_3\lambda)\in
E_{\delta=-y_1/y_2}(-\frac{x_1}{x_3},
-\frac{x_2}{x_3},\frac{1}{2x_3})\},&
\mathrm{if}\ \mathrm{Im}\,\frac{x_1}{x_3}< 0.
\end{array}
\right.
$$
Thus $E(x;y)$ consists of entire holomorphic
 functions $v(\lambda)$
obeying the resonance relation 
$v(-1+rx_3+sx_1)=
e^{2\pi is(2\sigma+x_2)/x_3}v(1+rx_3+sx_1)$, where $\sigma$ is the sign of $\mathrm{Im}\,
x_1/x_3$, 
and the bound
\begin{equation}\label{e-bound}
\left|v(\lambda)e^{\displaystyle
{-\frac{\pi i y_1\lambda^2}{2x_3y_2}}}
\right|
\leq C_1\exp\left(\pi
\frac{(\mathrm{Im}\,\frac\lambda{x_3})^2}
{\left|\mathrm{Im}\,\frac {x_1}{x_3}\right|}
+C_2|\lambda|\right),
\end{equation}
for some constants $C_1,C_2>0$.
Note that 
\[
E(x_1,x_2,x_3;y_1,y_2,y_3)
=E(-x_1,x_2,-x_3;-y_1,y_2,-y_3).
\]

For $\mathrm{Im} {x_1\over x_3} >0, \mathrm{Im} {x_2\over x_3} >0$,
introduce an integral operator
$U(x_1,x_2,x_3)$ acting on functions
in $E(x_2,x_1,x_3;y_2,y_1,y_3)$
as
\[
U(x_1,x_2,x_3)v(\lambda)=
{e^{{\pi i\over x_3}}\over 2\pi \sqrt{2ix_3}}
\int u\left(
\frac\lambda{x_3},
\frac\mu{x_3},
\frac{x_1}{x_3},
\frac{x_2}{x_3},
\frac{1}{2x_3}
\right)
Q\left(
\frac\mu{x_3},
\frac{x_2}{x_3},
\frac{1}{2x_3}
\right)
v(-\mu)d\mu.
\]
The integration path is given as follows. Suppose first
that $|x_2|,|x_3|>>1$. 
Then the integration is over a straight line which does not
intersect the segments joining $-1+rx_2+sx_3$ to
$1+rx_2+sx_3$ ($r,s\in\Z$) and such that the integrand
decays exponentially at infinity
along it. Such a path exists as the
integrand behaves as $\exp(-\mathrm{const}\mu^2)$ at infinity,
which decays for $\mu$ in a cone. 
The resonance relation implies that the integrand has opposite
residue at the pairs $\pm1+rx_2+sx_3$, so that,
up to sign, the integral is
independent  of the choice of path.
For $x_2,x_3$ general, the integral
is defined by analytic continuation.
The operator is defined up to multiplication by $\pm 1$ since one needs
to choose a square root and an orientation of the path.

We assume that 
\[
2\not\in\Z x_1+\Z x_2+\Z x_3.
\]
This ensures that $U(x_1,x_2,x_3)$ is well-defined. Indeed, the
integration kernel \eqref{eq-u} defining $U$ is regular at points obeying this
condition by the discussion of \ref{ss-onefactor} above. Moreover
this condition implies that
in the integral over
$\mu$ defining $U$, the integration contour is not pinched between
poles.

Proposition \ref{invert} says that
 $U(x_1,x_2,x_3)$
is an invertible linear map
\[
U(x_1,x_2,x_3):
E(x_2,x_1,x_3;y_2,y_1,y_3) \to
E(x_1,-x_2,x_3;y_1,-y_2,y_3),
\]
if $
\mathrm{Im}\,\frac{x_1}{x_3}>0,
\mathrm{Im}\,\frac{x_2}{x_3}>0,$

%

Introduce operators $\alpha(x_3)$, $\beta(x_1,x_2,x_3)$ where
\begin{eqnarray*}
\alpha(x_3)v(\lambda)&=& e^{-{i\pi\over x_3}
\left(\frac12\lambda^2-\frac13\right)}v(\lambda),
\\
\beta(x_1,x_2,x_3)v(\lambda)&=&
e^{-\frac{\pi i}{x_2x_3}\left( 
{x_1}\bigl(\frac1{2}\lambda^2-\frac13\bigr) + \lambda^2
-1\right)+\frac{2\pi i}{9x_1x_2x_3}}v(\lambda)
\qquad
\text{if}\, \text{Im}\, {x_2\over x_3}>0\,,
\\
\beta(x_1,x_2,x_3)v(\lambda)&=&
e^{-\frac{\pi i}{x_2x_3}\left( 
{x_1}\bigl(\frac1{2}\lambda^2-\frac13\bigr) - \lambda^2
+1\right)+\frac{2\pi i}{9x_1x_2x_3}}v(\lambda),
\qquad
\text{if}\, \text{Im}\, {x_2\over x_3}<0\,.
\end{eqnarray*}
The $\lambda$-independent terms in $\alpha$ and $\beta$ were
added to simplify the relations (1)-(4) below, which would 
otherwise only hold up to factors depending on $x$.

\begin{proposition}\label{p-alphabeta}
The operators $\alpha$ and $\beta$ are isomorphisms
\[
\alpha(x_3):
E(x_2,x_1,x_3;y_2,y_1,y_3)
\to
E(x_2,x_1+x_2,x_3;y_2-y_1,y_1,y_3),
\]
\[
\beta(x_1,x_2,x_3):
E(x_3,x_1,-x_2;y_3,y_1,-y_2)
\to
E(x_2,x_1,x_3,y_2,y_1,y_3),
\]
for all $(x,y)\in W$ for which they are defined.
\end{proposition}
\medskip
\noindent{\it Proof:} The statement for $\alpha$ is
easily checked inserting the definitions. The
fact that $\beta$ respects the
resonance relation is also easily checked.
We are left to prove that $v$ obeys the bound 
\eqref{e-bound}
if and only if $\beta v$ does. Let us prove the
only if part
in the case where $x_2/x_3$ has positive imaginary
part. The other cases are proved in the same way.
So we assume that 
$|\exp(\pi i y_3\lambda^2/2x_2y_1)v(\lambda)|
\leq C_1\exp(\pi\frac{(\mathrm{Im}\,\lambda/x_2)^2}
{\mathrm{Im}(-x_3/x_2)}+C_2|\lambda|)$. It follows
using the relation $\sum x_iy_i=0$ 
that $v'=\beta(x_1,x_2,x_3)v$ obeys the bound
 $|\exp(-\pi i y_2\lambda^2/2x_3y_1+\pi i\lambda^2/x_2x_3)
v'(\lambda)|
\leq C'_1\exp(\pi\frac{(\mathrm{Im}\,\lambda/x_2)^2}
{\mathrm{Im}(-x_3/x_2)}+C_2|\lambda|)$. 
The claim then follows from the identity
\[
\mathrm{Im}\,\frac{\lambda^2}{x_2x_3}
=
\frac{\left(\mathrm{Im}\,\frac\lambda{x_3}\right)^2}
{\mathrm{Im}\,\frac{x_2}{x_3}}
-
\frac{\left(\mathrm{Im}\,\frac\lambda{x_2}\right)^2}
{\mathrm{Im}\,-\frac{x_3}{x_2}}.
\]
\hfill$\square$

\medskip

Theorems \ref{q-heat-n}, \ref{t-mod}, 
\ref{thm-Fourier} take the form of identities:
\begin{enumerate}
\item 
 The  {\em q-heat equation},
\[
\alpha(x_3)\,U(x_1,x_1+x_2,x_3)\,\alpha(x_3)
\,U(x_1+x_2,x_2,x_3)\,\alpha(x_3)\,=\,  U(x_1,x_2,x_3)\,,
\]
holds on $E(x_2,x_1,x_3;y_2,y_1,y_3)$,
if $\mathrm{Im} {x_1\over x_3} >0, \mathrm{Im} {x_2\over x_3} >0$.
\item 
 The {\em first modular equation},
\begin{eqnarray*}
U(x_1,x_2,x_3)\,&&\beta (x_1,x_2,x_3)\,U(-x_3,-x_1,x_2)\,=
\\
\,&&\beta(-x_2,x_1,x_3)\,U(-x_3,x_2,x_1)\,\beta(-x_3,x_2,x_1)\,,
\end{eqnarray*}
holds on $E(-x_1,-x_3,x_2;-y_1,-y_3,y_2)$
if 
$\mathrm{Im} {x_1\over x_3} >0, \mathrm{Im} {x_2\over x_3} >0,
\mathrm{Im} {x_2\over x_1} >0$.
\item
  The {\em second modular equation,
\begin{eqnarray*}
U(-x_2,-x_3,x_1)\,&&\!\!\!\!\!\!\!\!\!
\beta (-x_2,-x_3,x_1)\,U(x_1,x_2,x_3)\,=
\\
\,&&\beta(x_3,x_2,-x_1)\,U(x_1,-x_3,x_2)\,\beta(x_1,x_3,-x_2)\,,
\end{eqnarray*}
 holds
on $E(x_2,x_1,x_3;y_2,y_1,y_3)$ if 
$\mathrm{Im} {x_1\over x_3} >0, \mathrm{Im} {x_2\over x_3} >0,
\mathrm{Im} {x_1\over x_2} >0$.}
\item
  The {\em inversion relation},
\begin{eqnarray*}
U(x_1,x_2,x_3)\,U(-x_2,-x_1,-x_3)=1\,,
\end{eqnarray*}
holds on $E(-x_1,-x_2,-x_3;-y_1,-y_2,-y_3)$
if 
$\mathrm{Im} {x_1\over x_3} >0, \mathrm{Im} {x_2\over x_3} >0$.
\end{enumerate}
Notice that each of the operators $U$ in these identities is defined
up to multiplication by $\pm 1$, so the right hand side of each of the
identities
is equal to the left hand side
up to multiplication by $\pm 1$.

\subsection{New range of parameters }\label{ss-newrange}

Extend the definition of the operator
 $U(x_1,x_2,x_3)$ from the domain
$\mathrm{Im} {x_1\over x_3} >0, \mathrm{Im} {x_2\over x_3} >0$ to the
domain $x_1/x_3, x_2/x_3 \in \C-\R$ by the formulas
\begin{eqnarray*}
U(x_1,x_2,x_3)&=& U(x_2,-x_1,x_3)^{-1} \qquad
\text{if} \,\qquad
\mathrm{Im} {x_1\over x_3} <0, \mathrm{Im} {x_2\over x_3} >0,
\\
U(x_1,x_2,x_3)&=& U(-x_2,x_1,x_3)^{-1} \qquad
\text{if} \,\qquad
\mathrm{Im} {x_1\over x_3} >0, \mathrm{Im} {x_2\over x_3} <0,
\\
U(x_1,x_2,x_3)&=& U(-x_1,-x_2,x_3) 
\qquad
\text{if} \qquad
\, \mathrm{Im} {x_1\over x_3} <0, \mathrm{Im} {x_2\over x_3} <0.
\end{eqnarray*}

\begin{thm}\label{t-hexagons} \ 
\begin{enumerate}
\item[(i)]
For any $(x;y)\in W$ such that $x_1/x_3, x_2/x_3 \in \C-\R$
and $2\not\in\Z x_1+\Z x_2+\Z x_3$,
the operator $U(x)$ defines an invertible linear map 
$E(x_2,x_1,x_3;y_2,y_1,y_3) \to
E(x_1,-x_2,x_3;y_1,-y_2,y_3)$.

\item[(ii)]

Moreover, the q-heat equation, first and second modular equations,
the inversion equation hold for $x_1,x_2,x_3$ such that $x_1/x_2,
x_1/x_3, x_2/x_3 \in \C-\R$ and $2\not\in\Z x_1+\Z x_2+\Z x_3$.

\end{enumerate}

\end{thm}

The six-term relations of 
Theorem \ref{t-hexagons} are the
commutativity of the following diagrams.

The {\em $q$-heat hexagon}:
\newcommand{\HEa} 
{E( x_1,\!-\!x_2, x_3; y_1,\!-\!y_2, y_3)}
\newcommand{\HEb}
{E( x_1,\!-\!x_1\!-\!x_2, x_3; y_1\!-\!y_2,\!-\!y_2, y_3)}
\newcommand{\HEc}
{E( x_2, x_1, x_3; y_2, y_1, y_3)}
\newcommand{\HEd}
{E( x_1\!+\!x_2, x_2, x_3; y_2,y_1\!-\!y_2, y_3)}
\newcommand{\HEe}
{E( x_2, x_1\!+\!x_2, x_3; y_2\!-\!y_1, y_1, y_3)}
\newcommand{\HEf} 
{E(x_1\!+\!x_2,\!-\!x_2, x_3;y_1,y_1\!-\!y_2,y_3)}
\begin{eqnarray*}
  & \HEc&  \\ 
&\alpha(x_3)\swarrow  
\hspace{3cm}\searrow 
U( x_1, x_2, x_3)  \\ 
& \HEe \hspace{3cm} \HEa\hspace*{1cm}& \\ 
&\left.U( x_1\!+\!x_2, x_2,x_3)
  \right\downarrow \hspace{5cm} 
\left\uparrow\alpha(x_3)\hspace*{2cm}\right.& \\ 
& \HEf \hspace{3cm}\HEb& \\
&\hspace*{1cm}\alpha(x_3)\searrow
  \hspace{3cm}\nearrow U( x_1, x_1\!+\!x_2, x_3) \\  &\HEd& 
\end{eqnarray*}

\medskip

The {\em first modular hexagon}:
\newcommand{\MEa} 
{E(-x_3, x_1, x_2;-y_3, y_1, y_2)}
\newcommand{\MEb}
{E( x_2, x_1, x_3; y_2, y_1, y_3)}
\newcommand{\MEc}
{E(-x_1,-x_3, x_2;-y_1,-y_3, y_2)}
\newcommand{\MEd}
{E( x_1,-x_2, x_3; y_1,-y_2, y_3)}
\newcommand{\MEe}
{E( x_2,-x_3, x_1; y_2,-y_3, y_1)}
\newcommand{\MEf} 
{E(-x_3,-x_2, x_1;-y_3,-y_2, y_1)}
\begin{eqnarray*}
  & \MEc&  \\ 
&\beta(-x_3, x_2, x_1)\swarrow  
\hspace{3cm}\searrow 
U(-x_3,-x_1, x_2)  \\ 
& \MEe \hspace{3cm} \MEa& \\ 
&\left.U(-x_3, x_2, x_1)
  \right\downarrow \hspace{5cm} 
\left\downarrow\beta( x_1, x_2, x_3)\right.& \\ 
& \MEf \hspace{3cm}\MEb& \\
&\beta(-x_2, x_1, x_3)\searrow
  \hspace{3cm}\swarrow U( x_1, x_2, x_3) \\  &\MEd& 
\end{eqnarray*}

\medskip

The {\em second modular hexagon}:
\newcommand{\MMEa} 
{E( x_1,-x_2, x_3; y_1,-y_2, y_3)}
\newcommand{\MMEb}
{E(-x_3,-x_2, x_1;-y_3,-y_2, y_1)}
\newcommand{\MMEc}
{E( x_2, x_1, x_3; y_2, y_1, y_3)}
\newcommand{\MMEd}
{E(-x_2, x_3, x_1;-y_2, y_3, y_1)}
\newcommand{\MMEe}
{E(-x_3, x_1, x_2;-y_3, y_1, y_2)}
\newcommand{\MMEf} 
{E( x_1, x_3, x_2; y_1, y_3, y_2)}
\begin{eqnarray*}
  & \MMEc&  \\ 
&\beta( x_1, x_3, -x_2)\swarrow  
\hspace{3cm}\searrow 
U( x_1, x_2, x_3)  \\ 
& \MMEe \hspace{3cm} \MMEa& \\ 
&\left.U( x_1,-x_3,  x_2)
  \right\downarrow \hspace{5cm} 
\left\downarrow\beta(-x_2, x_3, x_1)\right.& \\ 
& \MMEf \hspace{3cm}\MMEb& \\
&\beta( x_3, x_2,-x_1)\searrow
  \hspace{3cm}\swarrow U( -x_2,-x_3,  x_1) \\  &\MMEd& 
\end{eqnarray*}

The proof of Theorem \ref{t-hexagons} is done by
reduction to the case where the imaginary parts
of $x_1/x_3$ and $x_2/x_3$ are
positive. In fact it is 
straightforward to see that in all cases the identities for general imaginary parts can be rewritten
using the definitions as identities at some
other values of $x,y$ where the imaginary parts
are positive.

\subsection{Projectively flat $\SL(3,\Z)$ connections}\label{ss-profl}
We give here an interpretation of our relations in
terms of projectively flat discrete 
$\SL(3,\Z)$-connections.

We start by introducing the notion of discrete connections.
Let $X$ be a $G$-space with a fixed
presentation by generators
$e_j$, $j=1,\dots,k$ and relations $R_i=1$,
$i=1,\dots,r$.
A {\em discrete $G$-connection} on a 
complex vector bundle $\pi:F\to X$,
with fibers $F(x)=\pi^{-1}(x)$, $x\in X$,
assigns to each generator $e_j$
 a collection of linear isomorphisms $\phi_{e_j}(x):
F(e_j^{-1}x)\to F(x)$, $x\in X$. 
The {\em parallel translation} along an element $w$
of the free group $\mathrm{Free}_k$ 
generated by the $e_j$'s is the
collection of maps 
$\phi_w(x):F(\bar{w}^{-1}x)\to F(x)$ uniquely
defined by the properties
$\phi_{1}(x)=\mathrm{Id}_{F(x)}$, $\phi_{ww'}(x)=
\phi_w(x)\circ
\phi_{w'}(\bar w^{-1}x)$. Here $w\mapsto \bar w$
is the canonical projection $\mathrm{Free}_k\to G$.
The {\em curvature} of a discrete $G$-connection is 
the collection
of parallel translations $\phi_{R_i}(x)\in\mathrm{End}(F(x))$ along
the relations. A connection is called {\em projectively flat}
if $\phi_{R_i}(x)\in\C\,\mathrm{Id}$ 
for all $i$ and $x$.

Let $\pi:F\to Y$ be a vector bundle over a subset $Y$ of a $G$-space $X$.
Then a {\em discrete connection defined on $Y$} assigns to each generator $e_j$
 a collection of linear isomorphisms $\phi_{e_j}(x):
F(e_j^{-1}x)\to F(x)$, for all those
$x\in Y$ such that $e_{j}^{-1}x\in Y$. 
The parallel translation $\phi_w(x)$ is then defined on some subset
$Y_w\subset Y$. A discrete connection defined on $Y$ is projectively
flat if $\phi_{R_i}(x)\in\C\,\mathrm{Id}$ for all $i=1,\dots,r$ and
$x\in Y_{R_i}$.

In our case, $G=\SL(3,\Z)$ and for $X$ we take certain  orbits in $W$.

The group $\mathrm{SL}(3,\Z)$  is generated
by the elementary matrices 
$e_{ij}$, ($1\leq i, j\leq 3, i\neq j$) with ones in the 
diagonal and at $(i,j)$ and zeros everywhere else.
The relations
can be chosen  \cite{M} to be
\begin{eqnarray}\label{e-relations}
e_{ij}e_{kl}&=&e_{kl}e_{ij},\qquad i\neq l,\quad j\neq k,\notag\\
e_{ij}e_{jk}&=&e_{ik}e_{jk}e_{ij},\qquad i,j,k\; 
\mathrm{distinct},\\
(e_{13}\,e_{31}^{-1}e_{13})^4&=&1.\notag
\end{eqnarray}
An $\SL(3,\Z)$-orbit $X$ in $W$ is called {\em regular}
if for all $(x,y)\in X$\begin{enumerate}
\item $x_1/x_3\in\C-\R$ and $x_2/x_3\in\C-\R$ and
\item  $2\not\in\Z x_1+\Z x_2+\Z x_3$.
\end{enumerate}
Let $X$ be a regular orbit in $W$.
For  $(x,y)\in Y^\psi=\{(x,y)\in X\,|\,y_1\neq0\}$,
set
$$
F^\psi (x;y)=E(x_2,x_1,x_3; y_2,y_1,y_3).
$$
For  $(x,y)\in Y^\phi=\{(x,y)\in X\,|\,y_2\neq0\}$,
set
$$
F^\phi (x;y)=E(x_1,-x_2,x_3; y_1, -y_2,y_3).
$$

\begin{proposition}\label{p-iguanodon}  
Let $X\subset W$ be a regular orbit.
\begin{enumerate}
\item[(i)] The assignment $e_{ij}\mapsto \phi_{ij}(x)$ with
\begin{eqnarray*}
\phi_{12}(x)&=&\alpha(x_3)^{-1}U(x_1-x_2,x_1,x_3)^{-1}\alpha(x_3)^{-1},
\\
\phi_{13}(x)&=&\phi_{23}(x)=1,
\\
\phi_{21}(x)&=&\alpha(x_3)^{-1},
\\
\phi_{31}(x)&=&
\beta(-x_2,x_1-x_3,x_3),
\\
\phi_{32}(x)&=&\beta(x_3-x_2,x_3,-x_1)^{-1}
U(x_3-x_2,x_3,-x_1)^{-1}
\beta(-x_3,x_1,x_3-x_2)^{-1},
\end{eqnarray*}
defines a projectively flat 
$\SL(3,\Z)$-connection 
defined on $Y^{\phi}\subset X$.
\item[(ii)]
 The assignment $e_{ij}\mapsto \psi_{ij}(x)$ with
\begin{eqnarray*}
\psi_{12}(x)&=&\alpha(x_3),
\\
\psi_{13}(x)&=&\psi_{23}(x)=1,
\\
\psi_{21}(x)&=&\alpha(x_3)U(x_2,x_2-x_1,x_3)\alpha(x_3),
\\
\psi_{31}(x)&=&\beta(x_1-x_3,-x_3,x_2)^{-1}U(x_1-x_3,-x_3,x_2)^{-1}
\beta(x_3,x_2,x_3-x_1)^{-1},
\\
\psi_{32}(x)&=&
\beta(x_1,x_2-x_3,x_3),
\end{eqnarray*}
defines a projectively flat 
$\SL(3,\Z)$-connection 
defined on $Y^{\psi}\subset X$.
\item[(iii)] For all $i\neq j$ and all
$x\in Y^\phi\cap Y^\psi$ such that $e_{ij}^{-1}x\in Y^\phi\cap Y^\psi$, we have
\[
U(x)\psi_{ij}(x)=\phi_{ij}(x)U(e_{ij}^{-1}x).
\]
\end{enumerate}
\end{proposition}

\noindent{\it Proof:}
It first follows from Theorem \ref{t-hexagons} (i)
and Prop.\ \ref{p-alphabeta}
that 
$\phi_{ij}(x)$ is in all cases a well-defined isomorphism from
$F^\phi(e_{ij}^{-1}(x,y))$ to $F^\phi(x,y)$ and similarly for
$\psi$. The other claims are then simple consequences of 
the q-heat, modular and inversion relations (Theorem \ref{t-hexagons} (ii)), 
 and
the relations
\[
U(x_1,x_2+x_3,x_3)=U(x_1+x_3,x_2,x_3)=U(x_1,x_2,x_3),
\]
which follow from the fact that $u(\lambda,\mu,\tau,p,\eta)$ is
1-periodic in $\tau$ and $p$.

The easiest way to do the computations is to first check (iii), which
can be easily deduced from our three term relations. 
This identity implies that if the curvature associated to
a relation is scalar for one of the connections $\phi$ or $\psi$,
then it is scalar (and equal) also for the other connection.
Then the curvature can be computed using $\phi$
or $\psi$, whichever is simpler.
For example, to compute the curvature $C_{12}^{32}=\phi_{R}$ 
associated  to the relation
$R=e_{12}e_{32}e_{12}^{-1}e_{32}^{-1}$, it is better to
use the connection $\psi$. One gets 
\[
\psi_{12}(x)\psi_{32}(e_{12}^{-1}x)=C_{12}^{32}(x)
\psi_{32}(x)\psi_{12}(e_{32}^{-1}x),
\]
where
$$
C_{12}^{32}(x)=\exp\,\frac{2\pi i\, x_2}{9\,x_1x_3(x_1-x_2)(x_2-x_3)}\,.
$$
\hfill$\square$

In particular, if an orbit does not contain any point with
$y_1=0$ or $y_2=0$, which is true for generic orbits, we
have $Y^\phi=Y^\psi=X$ and the connections are defined 
everywhere. The trouble is that we do not know if the
spaces $F^\phi$, $F^\psi$ are nontrivial for these orbits.
By contrast, we have infinitely many examples of 
linearly independent functions
in $F^\psi(x,y)$ with $y_2=0$ or $F^\phi(x,y)$ with $y_1=0$.
Indeed, these spaces are isomorphic to spaces $E_\delta(\tau,-p,\eta)$
with $\delta=0$, which contain the functions
$u_\mu:\lambda \mapsto u(\lambda,\mu,\tau,p,\eta)$ with any fixed $\mu$,
see Prop.\ \ref{p-uinE0}. These functions are non-zero for generic $\mu$
since they are the meromorphic kernel of an invertible integral operator.
There are infinitely many linearly independent functions among them since
they behave differently under shifts by 1:
$u_\mu(\lambda+1)=
-e^{-\pi i\mu/2\eta}u_\mu(\lambda)$.
In this case we may construct a connection defined everywhere
by gluing the two partially defined connections: 

\begin{thm}\label{t-pf}
Let $X$ be a  regular orbit containing a point $(x,y)$ with
$y_1=0$. Let $Y^\psi$, $Y^\phi\subset X$ be the subsets as
above,
on which the connections $\phi$, $\psi$ are defined.
We have $Y^\phi\cup Y^\psi=X$.
Let $F$ be the vector bundle on $X$ obtained
from $F^\psi$ and $F^\phi$ by identifying the fibers
over $Y^\psi\cap Y^\phi$ via $U(x):F^\psi(x,y)\to F^\phi(x,y)$. 
Let a connection $\chi_{ij}(x)$ 
be defined as $\phi_{ij}(x)$, 
if both $(x,y)$ and $e_{ij}^{-1}(x,y)$ are in $Y^{\phi}$,
 and 
as $\psi_{ij}(x)$
if both $(x,y)$ and $e_{ij}^{-1}(x,y)$ are in $Y^{\psi}$. 
Then $\chi$ is a well-defined projectively
flat connection on $X$.
\end{thm}

\noindent{\it Proof:} Since any two points
 $(x,y),e_{ij}^{-1}(x,y)\in X$ related by a generator belong both
to $Y^{\phi}$ or both to $Y^{\psi}$, the connection $\chi_{ij}$  
is defined in all
cases. Moreover, by Prop.\ \ref{p-iguanodon} (iii), $\chi$
is well-defined on $F$, at least up to sign. 
It is easy  to check that that the curvatures along the relations
involve products of isomorphisms $\chi_{ij}(x)$
mapping between fibers at points $(x,y)$ 
which are all in $Y^{\phi}$ or all in $Y^\psi$.
Thus the claim that $\chi$ is a projectively flat connection
follows from the projective flatness of $\phi$ and $\psi$.
\hfill$\square$

\subsection{Formula for the curvature} We give here a 
formula for the curvature of the discrete connection $\phi$ (and
$\psi$). It is convenient to do a gauge transformation
$\bar\phi_{ij}(x)=\phi_{ij}(x)g_{ij}(x)$, with
\[
{g_{12}}(x)=\exp\left(\frac {2\pi i}{9\,
{x_{1}}\,({x_{1}} - {x_{2}}
)\,{x_{3}}}\right),\qquad
{g_{13}}(x)=\exp\left(\frac {2\pi i}{9\,{x_{1}}\,({x_{1
}} - {x_{3}})\,{x_{2}}}\right),
\]
\[
{g_{21}}(x)=\exp\left(  \frac {2\pi i}{9\,{x_{2}}\,({x_{2}} - {x_{
1}})\,{x_{3}}}\right),\qquad
{g_{23}}(x)=\exp\left(  \frac {2\pi i}{9\,{x_{2}}\,({x_{2}} - {x_{
3}})\,{x_{1}}}\right),
\]
\[
{g_{31}}(x)=1,\qquad
{g_{32}}(x)=\exp\left( -\, \frac {2\pi i}{9\,{x_{3}}\,({x_{3
}} - {x_{2}})\,{x_{1}}}\right).
\]
Then the curvature of the connection $\bar\phi$ has components
$\bar C_{ij}^{kl}(x)$, $\bar C(x)$ associated to the relations
\eqref{e-relations}. By definition, they are given by
\[
\bar\phi_{ij}(x)\bar\phi_{kl}(e_{ij}^{-1}x)=\bar C_{ij}^{kl}(x)
\bar\phi_{kl}(x)\bar\phi_{ij}(e_{kl}^{-1}x),\qquad i\neq l,\quad j\neq k,
\]
\[
\bar\phi_{ij}(x)\bar\phi_{jk}(e_{ij}^{-1}x)=\bar C_{ij}^{jk}(x)
\bar\phi_{ik}(x)\bar\phi_{jk}(e_{ik}^{-1}x)\bar\phi_{ij}(e_{jk}^{-1}e_{ik}^{-1}x),\qquad
i,j,k\; \text{distinct},
\]
\[
\bar\phi_s(x)\bar\phi_s(s^{-1}x)\bar\phi_s(s^{-2}x)\bar\phi_s(s^{-3}x)= \bar C(x),
\]
where $sx=(-x_3,x_2,x_1)$ and
\[
\bar\phi_s(x)=\bar\phi_{13}(x)\bar\phi_{31}(e_{31}e_{13}^{-1}x)^{-1}
\bar\phi_{13}(e_{13}^{-1}e_{31}x).
\]
We have
\[
\bar C_{12}^{32}(x)=\bar C_{32}^{12}(x)^{-1}=
\exp\left(\frac {2\pi i\,{x_{2}}}{3\,{x_{1}}\,{x
_{3}}\,({x_{2}} - {x_{3}})\,({x_{1}} - {x_{2}})}\right),
\]
\[
\bar C_{13}^{32}(x)=
\exp\left(\frac {2\pi i}{3\,{x_{3}}\,({x_{1}} - {x
_{2}})\,({x_{3}} - {x_{1}})}\right),\qquad
\bar C_{31}^{12}(x)=
\exp\left(\frac {2\pi i}{3\,{x_{1}}\,({x_{3}} - {x
_{2}})\,({x_{3}} - {x_{1}})}\right),
\]
and all other $\bar C_{ij}^{kl}(x)$ as well as $\bar C(x)$ are equal to 1.

\subsection{Comparison with the elliptic Gamma cocycle \cite{FV2}}
In \cite{FV2}, a nontrivial 2-cocycle of $\SL(3,\Z)$ with values in 
$\exp(2\pi i\Q(x_j/x_k)[z/x_3])$ was obtained by a similar construction
involving the elliptic gamma function instead of $U$. In the language of
discrete connections we use here, one defines an $\SL(3,\Z)$ connection
on the trivial line bundle over a dense subset of $\C^3$ by setting
\begin{eqnarray*}
\phi_{\Gamma,1,2}(x,z)&=&\Gamma\left(\frac {z-x_2}{x_3},
\frac{x_1-x_2}{x_3},-\,\frac{x_1}{x_3}\right)^{-1},\\
\phi_{\Gamma,3,2}(x,z)&=&\Gamma\left(\frac {z}{x_1},\frac{x_2-x_3}{x_1},
\frac{x_3}{x_1}\right),\\
\phi_{\Gamma,i,j}(x,z)&=&1,\qquad j\neq 2.
\end{eqnarray*}
For our purpose, $z$ may be considered here as a complex parameter on
which $\SL(3,\Z)$ acts trivially. The curvature of this connection
was computed in \cite{FV2}. It has the form $C_{\Gamma,ij}^{kl}=\exp(
\pi i L_{ij}^{kl}(z,x))$, $C_{\Gamma}(x)=1$, for some cubic polynomials
$L_{ij}^{kl}(x,z)$ in $z$ with coefficients in $\Q(x_1,x_2,x_3)$. The
relation to the curvature of $\bar\phi$ is
\[
\bar C_{ij}^{kl}(x)=\exp\left(-2\pi i\times\text{Coefficient of $z^3$ in $L_{ij}^{kl}(x,z)$}\right).
\]
The curvature of a projectively flat connection defines an extension
of the group and thus a characteristic class, see \cite{FV2}.
Conjecturally, $\bar C_{ij}^{kl}$ defines a nontrivial class in the
group cohomology $H^2(\SL(3,\Z),\exp(\C(x_1,x_2,x_3)))$. This means
(if the conjecture is true) that 
there is no gauge transformation given by exponentials of rational
functions that can make the curvature trivial.

\newpage
\appendix
\section{Theta functions and elliptic gamma functions}\label{App-A}

We summarize some formulae about theta functions, gamma functions
and phase functions, see \cite{FV2} for more details.
\subsection{The theta function}\label{A-theta}
Jacobi's first theta function is defined by the series
\[
\theta(z,\tau)=-\sum_{j\in\Z} e^{i\pi\tau(j+1/2)^2+2\pi
  i(j+1/2)(z+1/2)},\qquad z,\tau\in \C,\qquad \mathrm{Im}\,{}\tau>0.
\]
It is an entire holomorphic odd function such that
\begin{equation}\label{eq-frt}
  \theta(z+n+m\tau,\tau)=(-1)^{m+n}e^{-\pi im^2\tau-2\pi i
    mz}\theta(z,\tau), \qquad m,n\in\Z,
\end{equation}
and obeys the heat equation
\[
4\pi i\textstyle{\frac\partial{\partial\tau}}\theta(z,\tau)=
\theta''(z,\tau).
\]
Its transformation properties with respect to $\SL(2,\Z)$ are described
in terms of generators by the identities:
\[
\theta(-z,\tau)=-\theta(z,\tau),\qquad
\theta(z,\tau+1)=e^{\frac{i\pi}{4}}\theta(z,\tau), \qquad
\theta\left(\frac{z}{\tau}, -\,\frac{1}{\tau} \right)=i {\sqrt{-i\tau}}
e^{\frac{i\pi z^2}{\tau}} \theta(z,\tau).
\]
We choose the square root in the right half plane.

\subsection{Infinite products}
Let $x,q\in\C$ with $|q|<1$. The function
\[
(x;q)=\prod_{j=0}^\infty(1-xq^j)
\]
is a solution of the functional equation
\[
(qx;q)=\frac 1{1-x}(x;q).
\]
Let $x=e^{2\pi i z}$ and $q=e^{2\pi i \tau}$. Then
\begin{equation}\label{eq-thetafunction}
  \theta(z,\tau)=ie^{\pi i(\tau/4-z)}(x;q)(q/x;q)(q;q).
\end{equation}

We will also need the following variant of the theta function $\theta$:
\[
\theta_0(z,\tau)=(x;q)(q/x;q)=-i\frac{e^{\pi
    i(z-\tau/4)}}{(q;q)}\theta(z,\tau).\] This function obeys
\begin{eqnarray}\label{EGF1}
  \theta_0(z+1,\tau)&=& \theta_0(z,\tau),\nonumber \\ 
  \theta_0(z+\tau,\tau)&=& -e^{-2\pi iz}\theta_0(z,\tau), \\ 
  \theta_0(\tau-z,\tau)&=&\theta_0(z,\tau).\nonumber
\end{eqnarray}
Its modular properties 
are $\theta_0(z,\tau+1)=\theta_0(z,\tau)$, and if $z'=z/\tau$,
$\tau'=-1/\tau$,
\[
e^{\pi i(\tau/6-z)}\theta_0(z,\tau)= i\, e^{\pi
  i(-zz'+\tau'/6-z')}\theta_0(z',\tau').
\]
\subsection{Elliptic gamma functions}
Here we consider two parameters $\tau$ and $\sig$ in the upper half
plane, and set $q=e^{2\pi i \tau}$, $r=e^{2\pi i \sig}$, and consider
the function of $x=e^{2\pi i z}$,
\[
(x;q,r)=\prod_{j,k=0}^\infty(1-xq^jr^k)=(x;r,q).
\] It is a solution of the functional equations
\[
(qx;q,r)=\frac{(x;q,r)}{(x;r)}, \qquad (rx;q,r)=\frac{(x;q,r)}{(x;q)}.
\]
The {\em elliptic gamma function} \cite{R}, \cite{FV2}  is
\[
\Gamma(z,\tau,\sig)= \frac{(qr/x;q,r)}{(x;q,r)}= \Gamma(z,\sig,\tau).
\]
It obeys the identities
\begin{equation}\label{eq-Gamma9}
  \Gamma(z+1,\tau,\sig)=\Gamma(z,\tau,\sig),\qquad
  \Gamma(z+\sig,\tau,\sig)=\theta_0(z,\tau)\Gamma(z,\tau,\sig),
\end{equation}
and is normalized by $\Gamma((\tau+\sig+1)/2,\tau,\sig)=1$.  The zeros
of $(x;q,r)$ are at $x=q^{-j}r^{-k}$, $j,k=0,1,2,\dots$.  They are all
simple. Thus $\Gamma$ has only simple zeros and simple poles. The
zeros are at $ z=(j+1)\tau+(k+1)\sig+l, $ and the poles are at
$z=-j\tau-k\sig+l$. Here $j,k$ run over nonnegative integers and $l$
over all integers.

In fact, $\Gamma(z,\tau,\sig)$ is, up to normalization, the unique
1-periodic meromorphic solution of $u(z+\sig)=\theta_0(z,\tau)u(z)$
holomorphic in the upper half plane.
\subsection{Modular properties}
We consider the transformation properties of the elliptic gamma
function under modular transformations of $\sig$ and $\tau$. We have
the identities
\begin{eqnarray*}
  \Gamma(z,\tau,\sig)&=&\Gamma(z,\sig,\tau), \\ 
  \Gamma(z,\tau+1,\sig)&=&\Gamma(z,\tau,\sig), \\ 
  \Gamma(z,\tau+\sig,\sig) &=& \frac{\Gamma(z,\tau,\sig)}
  {\Gamma(z+\tau,\tau,\sig+\tau)}, \\ \Gamma(
  z/\sig,\tau/\sig,-1/\sig)&=& e^{i\pi Q(z;\tau,\sig)} {\Gamma(
    {({z-\sig})/\tau,-1/\tau,-\sig/\tau})} {\Gamma(z,\tau,\sig)} , \\ 
  Q(z;\tau,\sig)&=& \frac{z^3}{3\tau\sig} -
  \frac{\tau+\sig-1}{2\tau\sig} z^2 +
  \frac{\tau^2+\sig^2+3\tau\sig-3\tau-3\sig+1} {6\tau\sig} z \\ & & +
  \frac1{12} (\tau+\sig-1) (\tau^{-1}+\sig^{-1}-1).
\end{eqnarray*}
\subsection{Extending the range of parameters}\label{A-ext}
Since many operations we perform do not preserve the upper half plane,
it is important to extend the range of values $\tau$ and $\sig$ can
take.  We set
\[
(x;q^{-1})=\frac1{(qx;q)},\qquad (x;q^{-1},r)=\frac1{(qx;q,r)}, \qquad
(x;q,r^{-1})=\frac1{(rx;q,r)}.
\]
These formulae define an extension of the functions $(x;q)$, $(x;q,r)$
to meromorphic functions on $\{(x,q,r)| |q|\neq1\neq |r|\}$.  It is
clear that the functional relations
\[
(qx,q)=\frac1{1-x}(x,q),\qquad (qx;q,r)=\frac1{(x;r)}(x;q,r),
\]
still hold in this larger domain.  Correspondingly, we extend the
definition of $\theta_0$ and the elliptic gamma function by using the
same formulae in terms of the infinite products. We obtain:
\[\theta_0(z,-\tau)=\frac1{\theta_0(z+\tau,\tau)},\quad
\Gamma(z,-\tau,\sig)=\frac1{\Gamma(z+\tau,\tau,\sig)},\quad
\Gamma(z,\tau,-\sig)=\frac1{\Gamma(z+\sig,\tau,\sig)}.
\]
An easy check gives the following result:
\begin{proposition}
  All identities for $\Gamma$ and $\theta_0$ of the preceding
  subsections continue to hold for all $z,\tau,\sig$ such that
  $\tau,\sig\not\in\R$.
\end{proposition}

However, the statements about the position of zeros and poles are no
longer valid.

\subsection{The phase function}
We keep the notation of the previous subsection and introduce a new
variable $a$, and set $\alpha=e^{2\pi ia}$. The phase function is
\[
\Omega_a(z,\tau,\sig) =\frac{\Gamma(z+a,\tau,\sig)}
{\Gamma(z-a,\tau,\sig)} =\frac{(qr/x\alpha;q,r)(x/\alpha;q,r)}
{(x\alpha;q,r)(qr\alpha/x;q,r)}.
\]
We have
\begin{equation}\label{eq-zwz}
  \Omega_a(z+\sig,\tau,\sig) =\frac{\theta_0(z+a,\tau)}
  {\theta_0(z-a,\tau)} \Omega_a(z,\tau,\sig)
= e^{2\pi ia}\frac{\theta(z+a,\tau)}
  {\theta(z-a,\tau)} \Omega_a(z,\tau,\sig).
\end{equation}
The properties of this function follow from those of the gamma
function:
\begin{proposition}\label{p-1}
The function $\Omega_a(z,\tau,p)$ obeys the identities
\begin{equation}\label{eq-FROmega}
\Omega_a(z+p,\tau,p)=
e^{2\pi ia}\frac
{\theta(z+a,\tau)}
{\theta(z-a,\tau)}
\Omega_a(z,\tau,p),
\end{equation}
\begin{equation}\label{eq-FROmega2}
\Omega_a(z+\tau,\tau,p)=
e^{2\pi ia}\frac
{\theta(z+a,p)}
{\theta(z-a,p)}
\Omega_a(z,\tau,p),
\end{equation}
\begin{equation}\label{eq-FROmega3}
\Omega_a(z+1,\tau,p)=
\Omega_a(z,\tau,p),
\end{equation}
\begin{equation}\label{eq-Om2}
\Omega_a(z,\tau,p)=\Omega_a(z,p,\tau),
\end{equation}
\begin{equation}\label{eq-Om1}
\Omega_{a}(z,\tau,p)=
\Omega_{a}(z,\tau,\tau+p)
\Omega_{a}(z+p,\tau+p,p),
\end{equation}
 \begin{equation}\label{eq-Om1a}
\Omega_{a}(z,\tau,p)
=\Omega_{a}(z+\tau,\tau,\tau+p)
\Omega_{a}(z,\tau+p,p),
\end{equation}
\begin{equation}\label{eq-0m3}
\Omega_{a}(z,\tau+1,p)=
\Omega_{a}(z,\tau,p+1)
=\Omega_{a}(z,\tau,p),
\end{equation}
\begin{eqnarray}\label{eq-Om4}
\Omega_{a/\tau}\left(\frac z\tau,-\,\frac1\tau,\frac p\tau\right)
&=&
e^{\pi iS_a(z;\tau,p)}
\Omega_{a}(z,\tau,p)
\Omega_{a/p}(\frac {z-\tau}p,-\,\frac 1p,-\,\frac\tau p), \\
S_a(z;\tau,p)&=&\frac a{3\tau p}(6z^2\!-\!6(\tau\!+\!p\!-\!1)z\!+\!
2a^2\!+\!\tau^2\!+\!p^2\!+\!3\tau p\!-\!3\tau\!-\!3p\!+\!1),\nonumber
\end{eqnarray}
\begin{equation}\label{eq-Om5}
\Omega_a(z,\tau,p)\Omega_{-a}(z,\tau,p)=1,
\end{equation}
\begin{equation}\label{eq-Om6}
\Omega_a (-z,\tau,p)= \Omega_a (z,\tau,p)
\frac{\theta (z+a,\tau)\theta (z+a,p)}
{ \theta (z-a,\tau)\theta (z-a,p)}.
\end{equation}
\end{proposition}

\section{Modular properties of $R$-matrices and qKZB operators}\label{App-B}

Here we summarize some formulae giving the transformation properties
of the $R$-matrix. Fix $\Lambda_1,\Lambda_2\in\C$ and let
$R(z_1-z_2,\lambda,\tau,\eta)$ be the $R$-matrix of
$E_{\tau,\eta}(sl_2)$ associated to the evaluation Verma modules
$V_{\Lambda_1}(z_1)$, $V_{\Lambda_2}(z_2)$.  We have
\[
R(z_1-z_2,\lambda,\tau+1,\eta)=R(z_1-z_2,\lambda,\tau,\eta),
\]
and
\begin{eqnarray*}
  R\left( \frac z\tau, \frac \lambda\tau, -\,\frac 1\tau, \frac \eta\tau
  \right) &=& e^{(z_1-z_2) \frac{2\pi i\eta\Lambda_1\Lambda_2}{\tau} }
  A_1(\lambda-2\eta h^{(2)}) A_2(\lambda) \\ &&
\times\,  R(z_1-z_2,\lambda,\tau,\eta) A_1(\lambda)^{-1} A_2(\lambda-2\eta
  h^{(1)})^{-1}.
\end{eqnarray*}
Here $A_i(\lambda)=A(z_i,h^{(i)},\Lambda_i,\lambda,\eta)$ with
\[
A(z,h,\Lambda,\lambda,\eta)= \exp \frac{i\pi}\tau \left[
  z(h\lambda-\eta h^2) +\frac12(\Lambda-h) (\lambda+\eta\Lambda-\eta
  h) (\lambda-\eta\Lambda-\eta h) \right].
\]
These formulae can be deduced from the functional realization of
representations \cite{FTV1}: the $R$-matrix may be defined as the
unique linear map such that
\[
R(z_1-z_2,\lambda,\tau,\eta)\omega(z_1,z_2,\lambda,\tau,\eta)
=\omega^\vee(z_1,z_2,\lambda,\tau,\eta).
\]
The weight function $\omega=\omega_{ij}e^{\Lambda_1}_i\otimes
e^{\Lambda_2}_j$ and the mirror weight function $\omega^\vee$ take
values in the tensor product $V_{\Lambda_1}\otimes V_{\lambda_2}$ of
evaluation Verma modules. They are given explicitly in terms of ratios
of theta functions, see \cite{FTV1,FV1}.

{}From the above formulae, we deduce the transformation properties of
the qKZB operators. One finds
\begin{eqnarray*}
  K_i\left( \frac z\tau, \frac \lambda\tau, -\,\frac 1\tau, \frac p\tau,
    \frac \eta\tau\right) &\!\!\!=\!\!\! & e^{\frac{2\pi i\eta}{\tau}
    \sum_{j=1}^n \Lambda_i\Lambda_j (z_i-z_j+p\theta_{ij})}
  \prod_{j=1}^{n} A(z_j\!-\!p\theta_{ij},h^{(j)}, \Lambda_j,
  \lambda\!-\!2\eta\sum_{l=1}^{j-1}h^{(l)},\eta) \\ 
  &&\times\,K_i(z,\lambda,\tau,p,\eta)
  \prod_{j=1}^{n}A(z_j\!-\!p\theta_{ij},h^{(j)}, \Lambda_j,
  \lambda-2\eta\sum_{l=1}^{j-1}h^{(l)},\eta)^{-1},
\end{eqnarray*}
where we set $\theta_{ij}=1$ if $i>j$ and $\theta_{ij}=0$ if $i\leq
j$. In this somewhat simplified notation, the operators $A$ depend on
$\lambda$ and are to be viewed as multiplication operator by these
functions of $\lambda$; $K_i$ taken at $\lambda/\tau$ really means
that we conjugate $K_i$ by the operator of dilation of $\lambda$ by a
factor $1/\tau$.

\section{Estimates for theta functions}
\label{App-C}
\begin{lemma}\label{l-C1}
\begin{enumerate}
\item[(i)] For all $\tau$ in the upper half plane, there exists a
  constant $C_1(\tau)>0$ such that for all $\lambda\in\C$,
\[
|\theta(\lamdba,\tau)| \leq
C_1(\tau)\exp\left(\pi\frac{(\mathrm{Im}\,\lambda)^2}
  {\mathrm{Im}\,{}\tau} \right).
\]
\item[(ii)] For all $\epsilon>0$ and $\tau$ in the upper half plane,
  there exists a constant $C_2(\tau)>0$ such that for all $\lambda$
  such that $\mathrm{min}\{|\lambda-r-\tau s|,
  r,s\in\Z\}\geq\epsilon$,
\[
|\theta(\lamdba,\tau)| \geq C_2(\tau)\epsilon
\exp\left(\pi\frac{(\mathrm{Im}\,\lambda)^2}{\mathrm{Im}\,{}\tau}\right).
\]
\end{enumerate}
\end{lemma}

\medskip

\noindent{\it Proof:}
Let $S(\tau)$ be the compact set $S(\tau)= \{\lambda\in\C\,|\,
|\mathrm{Re}(\lambda)|\leq 1/2, |\mathrm{Im}\,{}\lambda| \leq
\mathrm{Im}\,{}\tau/2\}$. Then $\theta(\lambda,\tau)$ vanishes in
$S(\tau)$ only for $\lambda=0$.
 
We use the functional relation \Ref{eq-frt} of the theta function.  If
$\lambda=\lambda_0+m'+m\tau$, with $m',m\in\Z$ and $\lambda_0\in
  S(\tau)$ then
\begin{eqnarray*}
  |\theta(\lambda,\tau)|&=& e^{\pi m(2\mathrm{Im}\,{}\lambda_0+m
    \mathrm{Im}\,{}\tau)}|\theta(\lambda_0,\tau)| \\ &=& e^{\pi
    \frac{\mathrm{Im}\,{}\lambda-\mathrm{Im}\,{}\lambda_0}
    {\mathrm{Im}\,{}\tau}
    (\mathrm{Im}\,{}\lambda_0+\mathrm{Im}\,{}\lambda)}
  |\theta(\lambda_0,\tau)| \\ 
  &=&e^{\pi\frac{(\mathrm{Im}\,\lambda)^2}{\mathrm{Im}\,{}\tau}}
  e^{-\pi\frac{(\mathrm{Im}\,\lambda_0)^2}{\mathrm{Im}\,{}\tau}}
  |\theta(\lambda_0,\tau)|.
\end{eqnarray*}
Then the claim follows by setting
$C_1(\tau)=\max\{e^{-\pi\frac{(\mathrm{Im}\,\lambda_0)^2}{\mathrm{Im}\,{}\tau}}
|\theta(\lambda_0,\tau)|,\lambda_0\in S(\tau)\}$ and $C_2(\tau)=
\min\{e^{-\pi\frac{(\mathrm{Im}\,\lambda_0)^2}{\mathrm{Im}\,{}\tau}}
|\theta(\lambda_0,\tau)/\lambda_0|,\lambda_0\in S(\tau)\}$, so that
when $|\lambda_0|\geq \epsilon$, we have
\[
\theta(\lambda_0,\tau)\geq\theta(\lambda_0,\tau)\frac\epsilon{|\lambda_0|}\geq
C_2(\tau)\epsilon.
\]
\hfill$\square$

\end{document}